\documentclass[12pt,a4paper]{article}
\usepackage{amssymb,amsthm,amsmath}
\usepackage[english]{babel}
\usepackage{fullpage}
\numberwithin{equation}{section}
\title{Gauged Laplacians on quantum Hopf bundles}

\date{18 November 2008}

\author{~\\
\large{Giovanni Landi$^1$, Cesare Reina$^2$, Alessandro Zampini$^3$
\footnote{ Current address: Max Planck Institut f\"{u}r Mathematik, Vivatsgasse 7, D-53111 Bonn, Germany.  }
}
\\ [10pt]
\normalsize{$^1$ Dipartimento di Matematica e Informatica, Universit{\`a} di Trieste,} \\
\normalsize{Via  A. Valerio 12/1, I-34127, Trieste, Italy} \\
\normalsize{and INFN, Sezione di Trieste, Trieste, Italy} \\
\normalsize{ {\tt landi@univ.trieste.it} } \\ [10pt]
\normalsize{$^2$ Scuola Internazionale Superiore di Studi Avanzati,}\\
\normalsize{Via Beirut 2-4, I-34014, Trieste, Italy} \\
\normalsize{{\tt reina@sissa.it} } \\ [10pt]
\normalsize{$^3$ Institut f\"ur Angewandte Mathematik, Universit\"at Bonn,} \\
\normalsize{Wegelerstra\ss e 6, D-53115 Bonn, Germany} \\
\normalsize{ {\tt zampini@mpim-bonn.mpg.de}}
}

\newtheorem{theo}{Theorem}[section]
\newtheorem{lemm}[theo]{Lemma}
\newtheorem{prop}[theo]{Proposition}
\newtheorem{defi}[theo]{Definition}
\newtheorem{rema}[theo]{Remark}

\newcommand{\nn}{\nonumber}
\newcommand{\ce}{\mathcal{E}}
\newcommand{\dd}{{\rm d}}
\newcommand{\ca}{\mathcal{A}}

\newcommand{\ch}{\mathcal{H}}
\newcommand{\cl}{\mathcal{L}}
\newcommand{\cn}{\mathcal{N}}
\newcommand{\cp}{\mathcal{P}}
\newcommand{\cb}{\mathcal{B}}
\newcommand{\cq}{\mathcal{Q}}
\newcommand{\oca}[1]{\Omega^{#1}(\ca)}
\newcommand{\och}[1]{\Omega^{#1}(\ch)}
\newcommand{\cu}{\mathcal{U}}        %% an enveloping algebra
\newcommand{\SU}{\mathrm{SU}_q(2)}  %% quantum SU(2)
\newcommand{\ASU}{\ca(\mathrm{SU}_q(2))}  %% quantum SU(2)
\newcommand{\sq}{\mathrm{S}^2_{q}}  %% standard Podle\'s sphere
\newcommand{\Asq}{\ca(\mathrm{S}^2_{q})}  %% standard Podle\'s sphere
\newcommand{\su}{\cu_q(\mathrm{su}(2))}  %% quantum su(2)

\newcommand{\eps}{\varepsilon}      %% short for \varepsilon
\newcommand{\cop}{\Delta}           %% coproduct
\newcommand{\co}[2]{#1_{(#2)}}      %% coproduct factor : a_{(1)}
\newcommand{\hs}[2]{\left\langle #1,#2\right\rangle}  %% bilinear pairing

\newcommand{\ket}[1]{\left | #1 \right\rangle }
\newcommand{\bra}[1]{\left\langle #1 \right |}
\newcommand{\oh}{{\tfrac{1}{2}}}
\newcommand{\shalf}{{\scriptstyle\frac{1}{2}}} %% tiny fraction  1/2
\newcommand{\half}{{\mathchoice{\oh}{\oh}{\shalf}{\shalf}}} %% 1/2
\newcommand{\lt}{{\triangleright}}    %% a left action
\newcommand{\rt}{{\triangleleft}}
\newcommand{\IC}{{\mathbb C}} %% complex numbers
\newcommand{\IR}{{\mathbb R}} %% real numbers
\newcommand{\IN}{{\mathbb N}} %% natural numbers
\newcommand{\IZ}{{\mathbb Z}} %% integer numbers
\DeclareMathOperator{\Ad}{Ad}       %% ad map
\DeclareMathOperator{\id}{id}       %% identity map
\DeclareMathOperator{\Mat}{Mat}       %%  
\DeclareMathOperator{\End}{End}       %%  
\DeclareMathOperator{\SO}{SO}       %%  
\DeclareMathOperator{\U}{U}       %%  
\DeclareMathOperator{\tr}{tr}       %%  
\DeclareMathOperator{\qtr}{tr_q}       %%  
\DeclareMathOperator{\q-in}{q-index} 
\DeclareMathOperator{\chern}{ch}
\DeclareMathOperator{\q-in}{q-index}

\newcommand{\abs}[1]{\left|#1\right|}

\newcommand{\figureheight}{8cm}
\newcommand{\putfig}[2]{\begin{figure}[htp]
        \special{isoscale c:/itex/texfig/#1.wmf, \the\hsize \figureheight}
        \vspace{\figureheight}
        \caption{#2}\label{fig:#1}
        \end{figure}}
\newcommand{\pictureheight}{4cm}
\newcommand{\putpicture}[2]{\begin{figure}[htp]
        \special{isoscale c:/itex/texfig/#1.wmf, \the\hsize \pictureheight}
        \vspace{\pictureheight}
        \caption{#2}\label{fig:#1}
        \end{figure}}

\newcommand{\beqa}{\begin{eqnarray}}
\newcommand{\eeqa}{\end{eqnarray}}
\newcommand{\beq}{\begin{equation}}
\newcommand{\eeq}{\end{equation}}

\newcommand{\del}{\partial}

\newcommand{\mn}{\abs{n}}
\newcommand{\qpp}{\mathfrak{p}^{\left(n\right)}}
\newcommand{\qpn}{\check{\mathfrak{p}}^{\left(n\right)}}
\newcommand{\bz}{B_{0}}
\newcommand{\bp}{B_{+}}
\newcommand{\bm}{B_{-}}
\newcommand{\delb}{\bar{\del}}

\begin{document}
\maketitle

\thispagestyle{empty}

\maketitle

\begin{abstract}
\noindent 
We study gauged Laplacian operators on line bundles on a quantum 2-dimensional sphere. Symmetry under the (co)-action of a quantum group allows for their complete diagonalization.  These operators describe `excitations moving on the quantum sphere' in the field of a magnetic monopole. The energies are not invariant under the exchange monopole/antimonopole, that is under inverting the direction of the magnetic field. There are potential applications to models of quantum Hall effect. 
%with quantum-group symmetries. 
\end{abstract}

\newpage
\section{Introduction}
The motivation for the present paper is two-fold consisting in generalizing to noncommutative manifolds and bundles over them important mathematical and physical models. 

On the mathematical side we start from the following classical construction. 
Let $(M,g)$ be a compact Riemannian manifold and 
$P\to M$ a principal bundle with compact structure Lie group $G$.  With  
$(\rho, V)$  a representation of $G$, there is the well known identification of sections of the associated vector bundle $E = P\times_G V$ on $M$ with equivariant maps from $P$ to $V$, $\Gamma(M,E) \simeq  C^\infty(P,V)_G \subset C^\infty(P) \otimes V$. Given a connection on $P$ one has a covariant derivative $\nabla$ on $\Gamma(M,E)$, that is $\nabla: \Gamma(M,E) \to \Gamma(M,E) \otimes_{C^\infty(M)}\Omega^1(M)$. 
The Laplacian operator, $\Delta^E := -(\nabla \nabla^* + \nabla^*\nabla)$, 
where the dual $\nabla^*$ is defined with the metric $g$, is a map from $\Gamma(M,E)$ to itself. This operator is related to the Laplacian operator  
on the total space $P$, $\Delta^P = -(\dd \dd^* + \dd^* \dd)$ acting on $C^\infty(P)$, as
$$
\Delta^E = \big (\Delta^P \otimes 1 + 1 \otimes C_G \big )_{\big | C^\infty(P,V)_G}
$$ 
(see e.g. \cite[Prop.~5.6]{bgv}). 
Here $C_G = \sum_a \rho(\xi_a)^2 \in \End(E)$ is the quadratic Casimir operator of $G$, with $\{\xi_a, a=1, \dots, \dim G\}$ an orthonormal basis of the Lie algebra of $G$. The metric on $P$ used for the dual $\dd^*$  is the canonical one obtained from the metrics on $M$ and on the Lie algebra of $G$. 
We refer to $\Delta^E$ as the gauged Laplacian of $M$ in the sense that the Laplacian $\Delta^M$ of $M$ gets replaced by $\Delta^E$ when  the exterior derivative $\dd$ is replaced by the connection $\nabla$. The above formula has several important applications, notably the study of the heat kernel expansion and index theorems on principal bundles. 

Our mathematical motivation is then to look for formul{\ae} like the above relating Laplacian operators on noncommutative bundles and study their possible use.

For the physical motivation the starting observation is the fact that the Laughlin wave functions \cite{L83} for the fractional quantum Hall effect (on the plane) are not translationally invariant. 
This problem was overcome in \cite{H83} with a model on a sphere with a magnetic monopole at the origin. 
The full Euclidean group of symmetries of the plane is recovered from the rotation group $\SO(3)$ of symmetries of the sphere. In fact, one is really dealing with the Hopf fibration of the sphere $\mathrm{S}^3=\mathrm{SU(2)}$ over the sphere $\mathrm{S}^2$ with $\U(1)$ as gauge (or structure) group and needs to diagonalize the gauged Laplacian of $\mathrm{S}^2$ gauged with the monopole connection. For this bundle one has that $\Delta^P=\Delta^{\mathrm{SU(2)}}=C_{\mathrm{SU(2)}}$ and by the formula above the diagonalization is straightforward. 

It is then natural to seek models of Hall effect on noncommutative spaces. These will enjoy symmetry for quantum groups and lead  to potentially interesting 
physical models. In the present paper we study a model of `excitations moving on a quantum 2-sphere' and in the field of a magnetic monopole. The most striking fact is that the energies of the corresponding gauged Laplacian operator are not invariant under the exchange monopole/antimonopole,
that is under inversion of the direction of the magnetic field, a manifestation of the phenomenon that `quantization removes degeneracy'

More specifically, the paper is organized as follows.
In Sect.~\ref{se:qhb} we describe  a quantum principal $\U(1)$-bundle over a quantum sphere $\sq$ having as total space the manifold of the quantum group $\SU$; there are natural associated line bundles classified by the winding number $n\in\IZ$. We describe in Sect.~\ref{se:wn} the computation of the winding number of the bundles as well as a `twisted' version of it.  On the bundle one introduces differential calculi -- recalled in Sect.~\ref{se:cotqpb} -- that lead to the notion of gauge connection, given in Sect.~\ref{se:con}. The gauged Laplacian operator on sections of associated bundles in introduced in Sect.~\ref{se:glo}. As mentioned, this operator describes `excitations moving on the quantum sphere' and in the field of a magnetic monopole. There is symmetry under the action of the quantum universal enveloping algebra $\su$ -- the Hopf algebra dual to the quantum group  $\SU$ -- that allows for a complete diagonalization of the gauged Laplacian. 
The energies of the Laplacian depend explicitly on the deformation parameter, an improvement with respect to a similar model of excitations in the field of instantons on a noncommutative 4-sphere \cite{La06} and, as mentioned, they are not invariant under the exchange monopole/antimonopole. The relation of the gauged Laplacian $\Box_{\nabla}$ with the quadratic Casimir operator $C_{q}$ of $\su$ is more involved than the classical one:
$$
q K^2 \Box_{\nabla} = C_{q} + \tfrac{1}{4} - \half 
\left(\frac{qK^{2}-2+q^{-1}K^{-2}}{(q-q^{-1})^{2}} +
\frac{q^{-1}K^{2}-2+qK^{-2}}{(q-q^{-1})^{2}}\right) 
\, ,
$$
with $K$ the group like element of $\su$, the `generator of the structure group' $\U(1)$. 
In order to have a reasonable self-contained paper while making an effort to lighten it, we give in App.~\ref{se:pre} some general constructions on differential calculi and quantum principal bundles with connections, and relegate to App.~\ref{se:cqpb} some of the computations relevant to the $\U(1)$-bundle over the quantum sphere $\sq$ of Sect.~\ref{se:qhb}.  We leave to future work a more detailed study of potential applications to models of quantum Hall effect (with quantum group symmetries) as well as to heat kernel methods for noncommutative index theorems.

\section{The quantum Hopf bundle}\label{se:qhb}

The quantum principal bundle we need is the well know $\U(1)$-bundle  over the standard Podle\'s sphere $\sq$ and whose total space is the manifold of the quantum group $\SU$. This bundle is an example of a quantum homogeneous space \cite{BM93}. We recall in App.~\ref{se:pre} the general framework of quantum principal bundles with nonuniversal calculi and endowed with a (gauge) connection. In this section we start with the algebras of the total and base space of the bundle and then we proceed to describe all bundles associated with the action of the group $\U(1)$. 
We also describe at length the related Peter-Weyl decomposition for the algebra $\ASU$ that we shall need later on in Sect.~\ref{se:glo}. 

%In this section we present the quantum Hopf bundle as an example of a quantum principal bundle. In particular, we want to show the well-known fact that such a quantum bundle is an example of a quantum homogeneous bundle \cite{brz94,BM93}, whose total space is the manifold of the quantum group $\SU$  over the standard Podle\'s
%sphere $\sq$ with $\U(1)$ as structure group.

\subsection{The algebras}\label{qdct}

The coordinate algebra $\ASU$ of the quantum group
$\SU$ is the $*$-algebra generated by $a$ and~$c$, with relations,
\begin{align}\label{derel}
& ac=qca\,\,\,\,\,ac^*=qc^*a\,\,\,\,\,cc^*=c^*c , \nn \\
& a^*a+c^*c=aa^*+q^{2}cc^*=1 .
\end{align}
The deformation parameter $q\in\IR$ is taken in the interval 
$0<q<1$, since for $q>1$ one gets isomorphic algebras; at $q=1$ one recovers the commutative coordinate algebra on the manifold $\mathrm{SU(2)}$.
The Hopf algebra structure for $\ASU$ is given by the coproduct:
$$
\Delta\,\left[
\begin{array}{cc} a & -qc^* \\ c & a^*
\end{array}\right]=\left[
\begin{array}{cc} a & -qc^* \\ c & a^*
\end{array}\right]\otimes\left[
\begin{array}{cc} a & -qc^* \\ c & a^*
\end{array}\right] ,
$$
antipode:
$$
S\,\left[
\begin{array}{cc} a & -qc^* \\ c & a^*
\end{array}\right]=\left[
\begin{array}{cc} a^* & c^* \\ -qc & a
\end{array}\right],
$$
and counit:
$$\epsilon \left[
\begin{array}{cc} a & -qc^* \\ c & a^*
\end{array}\right]=\left[
\begin{array}{cc} 1 & 0 \\ 0 & 1
\end{array}\right].
$$

\bigskip
The quantum universal enveloping algebra $\su$ is the Hopf $*$-algebra
%over $\complex$
generated as an algebra by four elements $K,K^{-1},E,F$ with $K K^{-1}=1$ and  subject to
relations: 
\beq \label{relsu} 
K^{\pm}E=q^{\pm}EK^{\pm}, \qquad 
K^{\pm}F=q^{\mp}FK^{\pm}, \qquad
[E,F] =\frac{K^{2}-K^{-2}}{q-q^{-1}}. 
\eeq 
The $*$-structure is simply
$$
K^*=K, \qquad E^*=F, \qquad F^*=E, 
$$
and the Hopf algebra structure is provided  by coproduct $\Delta$, antipode $S$, counit $\epsilon$: 
$$
\begin{array}{cc}
\Delta(K^{\pm}) =K^{\pm}\otimes K^{\pm},
\qquad \Delta(E) =E\otimes K+K^{-1}\otimes E, \qquad \Delta(F)
=F\otimes K+K^{-1}\otimes F , 
\\~ \\
S(K) =K^{-1},
\qquad S(E) =-qE, \qquad S(F) =-q^{-1}F , 
\\~\\
\epsilon(K)=1, \qquad \epsilon(E)=\epsilon(F)=0 .
\end{array}
$$
{}From the
relations \eqref{relsu}, the quadratic quantum Casimir element:
\beq\label{cas}
C_{q}\,:=\,\frac{qK^{2}-2+q^{-1}K^{-2}}{(q-q^{-1})^{2}}+FE-\tfrac{1}{4}
\eeq
generates the centre of $\su$. In the `classical limit $q \to 1$' the operator $C_{q}$ goes to the Casimir element $C_{\mathrm{SU(2)}}=H^2+\half(EF+FE)$, as can be seen by setting $K=q^H$, expanding in the parameter $\hbar=:\log q$ and truncating at the 0-th order in $\hbar$. In the limit the elements $H,E,F$ are the generators of the Lie algebra $\mathrm{su}(2)$.

There is a bilinear pairing between $\su$ and $\ASU$, given on
generators by
\begin{align*}
&\langle K,a\rangle=q^{-1/2}, \quad \langle K^{-1},a\rangle=q^{1/2}, \quad
\langle K,a^*\rangle=q^{1/2}, \quad \langle K^{-1},a^*\rangle=q^{-1/2}, \nn\\
&\langle E,c\rangle=1, \quad \langle F,c^*\rangle=-q^{-1},
\end{align*}
and all other couples of generators pairing to~0. One regards $\su$ as a subspace of the linear dual of~$\ASU$ via this
pairing. There are \cite{wor87} canonical left and right $\su$-module algebra
structures on~$\ASU$  such that
$$
\hs{g}{h \lt x} := \hs{gh}{x},  \quad  \hs{g}{x \rt  h} := \hs{hg}{x},
\qquad \forall\, g,h \in \su,\ x \in \ASU.
$$
They are given by $h \lt x := \hs{(\id \otimes h)}{\cop x}$ and
$x \rt  h := \hs{(h \otimes \id)}{\cop x}$, or equivalently, 
$$
h \lt x := \co{x}{1} \,\hs{h}{\co{x}{2}}, \qquad
x \rt  h := \hs{h}{\co{x}{1}}\, \co{x}{2},
$$
in the Sweedler notation.
These right and left actions are mutually commuting:
$$
(h \lt a) \rt  g = \left(\co{a}{1} \,\hs{h}{\co{a}{2}}\right) \rt  g 
= \hs{g}{\co{a}{1}} \,\co{a}{2}\, \hs{h}{\co{a}{3}} 
= h \lt \left(\hs{g}{\co{a}{1}}\, \co{a}{2}\right) = h \lt (a \rt  g),
$$
and since the pairing satisfies 
$$
\hs{(Sh)^*}{x} = \overline{\hs{h}{x^*}},
\qquad \forall\, h \in \su,\ x \in \ASU ,
$$
the $*$-structure is compatible with both actions:
$$
h \lt x^* = ((Sh)^* \lt x)^*,  \quad  x^* \rt  h = (x \rt  (Sh)^*)^*,
\qquad \forall\,  h \in \su, \ x \in \ASU.
$$
We list both actions on powers of generators. Here and below we shall use the `$q$-number',
\begin{equation}
[x] = [x]_q := \frac{q^x - q^{-x}}{q - q^{-1}} ,
\label{eq:q-integer}
\end{equation}
defined for $q \neq 1$ and any $x \in \IR$. First the
left action; for $s=1,2,\ldots $:
\begin{align}\label{lact}
& K^{\pm}\triangleright a^{s} =q^{\mp\frac{s}{2}}a^{s}, \quad
K^{\pm}\triangleright a^{* s} =q^{\pm\frac{s}{2}}a^{* s}, \quad
K^{\pm}\triangleright c^{s} =q^{\mp\frac{s}{2}}c^{s}, \quad
K^{\pm}\triangleright c^{* s} =q^{\pm\frac{s}{2}}c^{* s} ;  \\
& \nn\\
& F\triangleright a^{s} =0, \quad  F\triangleright
a^{*s} =q^{(1-s)/2} [s] c a^{* s-1}, \quad
F\triangleright c^{s} =0, \quad
F\triangleright c^{*s} =-q^{-(1+s)/2} [s] a c^{*s-1} ; \nn \\
& \nn\\
& E\triangleright a^{s} =-q^{(3-s)/2} [s] a^{s-1} c^{*} , \quad
E\triangleright a^{* s} =0 , \quad
E\triangleright c^{s} =q^{(1-s)/2} [s]  c^{s-1} a^*, \quad
E\triangleright c^{* s} =0 . \nn
\end{align}
Then the right one:
\begin{align}\label{ract}
& a^{s}\triangleleft K^{\pm} =q^{\mp\frac{s}{2}}a^{s}, \quad
a^{* s}\triangleleft K^{\pm} =q^{\pm\frac{s}{2}}a^{* s}, \quad
c^{s}\triangleleft K^{\pm} =q^{\pm\frac{s}{2}}c^{s}, \quad
c^{* s}\triangleleft K^{\pm} =q^{\mp\frac{s}{2}}c^{* s} ; \\
& \nn\\
& a^{s}\triangleleft F =q^{(s-1)/2} [s] c a^{s-1}, \quad
a^{* s}\triangleleft F =0 , \quad
c^{s}\triangleleft F =0 , \quad
c^{* s}\triangleleft F =-q^{-(s-3)/2} [s] a^{*}c^{*s-1} ; \nn \\
& \nn\\
& a^{s}\triangleleft E =0, \quad
a^{*s}\triangleleft E =-q^{(3-s)/2} [s] c^{*}a^{*s-1}, \quad
c^{s}\triangleleft E =q^{(s-1)/2} [s] c^{s-1} a, \quad
c^{* s}\triangleleft E =0 . \nn
\end{align}

\bigskip
We also need to recall (see e.g. \cite[Prop.~3.2.6]{maj95}) that the irreducible finite dimensional
$*$-representations $\sigma_J$ of $\su$ are labelled by nonnegative
half-integers $J\in \half \IN$ (the spin); they are given by
\begin{align}
\sigma_J(K)\,\ket{J,m} &= q^m \,\ket{J,m},
\nn \\
\sigma_J(E)\,\ket{J,m} &= \sqrt{[J-m][J+m+1]} \,\ket{J,m+1},
\label{eq:uqsu2-repns} \\
\sigma_J(F)\,\ket{J,m} &= \sqrt{[J-m+1][J+m]} \,\ket{J,m-1},
\nn
\end{align}
where the vectors $\ket{J,m}$, for $m = J, J-1,\dots, -J+1, -J$, form
an orthonormal basis for the $(2J+1)$-dimensional, 
irreducible $\su$-module $V_J$, and the brackets
denote the $q$-number as in~\eqref{eq:q-integer}. Moreover, $\sigma_J$
is a $*$-representation of $\su$, with respect to the hermitian
scalar product on $V_J$ for which the vectors $\ket{J,m}$ are
orthonormal. In each representation $V_J$, the Casimir \eqref{cas}
is a multiple of the identity with constant given by:
\beq
C_{q}^{(J)}= [J+\half]^2- \tfrac{1}{4}.
\label{quca}
\eeq

\bigskip
Denote $\ca(\U(1)):=\IC[z,z^*] \big/ \!\!<zz^* -1>$; the map:
\beq  \label{qprp}
\pi: \ASU \, \to\,\ca(\U(1)),
\qquad \pi\,\left[
\begin{array}{cc} a & -qc^* \\ c & a^*
\end{array}\right]=
\left[
\begin{array}{cc} z & 0 \\ 0 & z^*
\end{array}\right]
\eeq 
is a surjective Hopf $*$-algebra homomorphism, so that $\ca(\U(1))$
becomes a quantum subgroup of $SU_{q}(2)$ with a right coaction,
\beq 
\Delta_{R}:= (\id\otimes\pi) \circ \Delta \, : \, \ASU \,\mapsto\,\ASU \otimes
\ca(\U(1)) . \label{cancoa} 
\eeq 
The coinvariant elements for this
coaction, that is
elements $b\in\ASU$ for which $\Delta_{R}(b)=b\otimes 1$,  
form a subalgebra of $\ASU$ which is the coordinate algebra
$\Asq$ of the standard Podle\'s sphere $\sq$. It is straightforward
to establish that 
$$
\Delta_{R}(a)=a\otimes z , \quad
\Delta_{R}(a^*)=a^*\otimes z^{*} , \quad \Delta_{R}(c)=c\otimes z ,
\quad \Delta_{R}(c^*)=c^*\otimes z^{*} . 
$$
As a set of generators (different from the original ones in \cite{Po87})
for $\Asq$ we take 
\beq \label{podgens}
B_{-} := - ac^* , \qquad 
B_{+} :=q ca^* , \qquad 
B_{0} :=  \frac{q^{2}}{1+q^{2}} - q^{2} cc^*
,\eeq 
for which one finds relations:
\begin{align*}
& B_{-}B_{0} = q^{2} B_{0}B_{-}, \qquad B_{0}B_{+} = q^{2} B_{+}B_{0}, \\
& B_{+}B_{-}=q \left[q^{-2} B_{0} - (1+q^{2})^{-1} \right] 
\left[q^{-2} B_{0} + (1+q^{-2})^{-1} \right] , \\
& B_{-}B_{+}=q \left[B_{0} + (1+q^{2})^{-1} \right] 
\left[B_{0} - (1+q^{-2})^{-1} \right],
\end{align*}
and $*$-structure:
\[
(B_{0})^*=B_{0}, \qquad (B_{+})^*= - q B_{-} .
\]
The sphere $\sq$ is a
quantum homogeneous space of $\SU$ and the coproduct of $\ASU$
restricts to a left coaction of $\ASU$ on $\Asq$ which, on
generators reads:
\begin{align*}
\Delta(B_{-})&=a^{2}\otimes B_{-}-(1+q^{-2})B_{-}\otimes B_{0}+
c^{* 2}\otimes B_{+} , \nn\\
\Delta(B_{0})&= q\, ac\otimes B_{-}+(1+q^{-2})
B_{0}\otimes B_{0}- c^*a^* \otimes B_{+} , \nn\\
\Delta(B_{+})&=q^{2}\, c^{2}\otimes B_{-}+(1+q^{-2})B_{+}\otimes
B_{0}+a^{* 2}\otimes B_{+} .
\end{align*}
The algebra inclusion $\Asq\hookrightarrow\ASU$ is a quantum principal bundle and can be endowed with compatible calculi \cite{BM93}, a construction that we shall illustrate below.

\subsection{The associated line bundles}\label{se:avb}

The left action of the group-like element $K$ on $\ASU$ allows one \cite[eq.~(1.10)]{maetal} to give a
vector basis decomposition $\ASU=\oplus_{n\in\IZ} \cl_n$, where, 
\beq\label{libu} 
\cl_n := \{x \in \ASU ~:~ K \lt x = q^{n/2} x \} .
\eeq 
In particular $\Asq = \cl_0$.  Also, $\cl_n^* \subset \cl_{-n}$ and $\cl_n\cl_m \subset \cl_{n+m}$. 
Each $\cl_n$ is clearly a bimodule over $\Asq$.
It was shown in \cite[Prop.~6.4]{SWPod} that
each $\cl_n$ is isomorphic to a projective left $\Asq$-module of rank 1.
These projective left $\Asq$-modules
give modules of equivariant maps or of sections of line bundles over the quantum sphere $\sq$ with winding numbers (monopole charge) $-n$ (see Sect.~\ref{se:wn} below).
The corresponding projections (cf. \cite{BM98,HM98}) can be written explicitly.
They are elements $\qpp$ and $\qpn$ in $\Mat_{\mn+1}(\Asq)$ (for $n\geq0$ and $n\leq0$ respectively) whose explicit form we give below. 
For $n\geq0$, 
\begin{align}\label{qpro}
& \qpp = \ket{\Psi^{(n)}} \bra{\Psi^{(n)}} \, , \qquad \qquad 
\ket{\Psi^{(n)}}_{\mu} = \sqrt{\beta_{n,\mu}} ~ c^{* \mu}a^{* n-\mu} \, ,
\\
&  
\mathrm{with} \qquad \beta_{n,\mu}=q^{2\mu}\prod\nolimits_{j=0}^{\mu-1}\left(\frac{1-q^{-2\left(n-j\right)}}{1-q^{-2\left(j+1\right)}}\right),
\quad \mu = 0, \dots,  n \, . \nn
\end{align}
Here and below we use the convention that $\prod_{j=0}^{-1}(\cdot)=1$.
The entries of $\qpp$ are:
\[
\qpp_{\mu\nu}=\sqrt{\beta_{n,\mu}\beta_{n,\nu}}\,
c^{*\mu}a^{*n-\mu}a^{n-\nu}c^{\nu} \in \Asq.
\]
For $n\leq0$ one has instead,
\begin{align}\label{qpron}
& \qpn = \ket{\check{\Psi}^{(n)}} \bra{\check{\Psi}^{(n)}}  \, , \qquad \qquad
\ket{\check{\Psi}^{(n)}}_{\mu} =  \sqrt{\alpha_{n,\mu}}  ~ c^{\mn-\mu}a^{\mu} \, , 
 \\
& 
\mathrm{with} \qquad \alpha_{n,\mu}=\prod\nolimits_{j=0}^{\mn-\mu-1}\left(\frac{1-q^{2\left(\mn-j\right)}}
{1-q^{2\left(j+1\right)}}\right),  \quad \mu = 0, \dots,  \mn \, . \nn
 \end{align}
The entries of $\qpn$ are:
\[
\qpn_{\mu\nu}= \sqrt{\alpha_{n,\mu}\alpha_{n,\nu}}\,c^{\mn-\mu}a^{\mu}
a^{*\nu}c^{*\mn-\nu} \in \Asq .
\]
Both $\qpp$ and $\qpn$ are self-adjoint by construction.
Using the commutation relations \eqref{derel} of $\ASU$ and the
explicit form of the coefficients \eqref{qpro} and \eqref{qpron} a
long but straightforward computation shows that: 
\beq\label{id1}
\hs{\Psi^{(n)}}{\Psi^{(n)}} = 
\sum\nolimits_{\mu=0}^{n} 
 \beta_{n,\mu}\, a^{n-\mu} c^{\mu} c^{* \mu} a^{* n-\mu} =
 (aa^{*}+q^{2}cc^{*})^{n} = 1 ,
\eeq
and analogously
\beq\label{id2}
 \hs{\check{\Psi}^{(n)}}{\check{\Psi}^{(n)}} =
\sum\nolimits_{\mu=0}^{\mn} 
 \alpha_{n,\mu}\, a^{* \mu} c^{* \mn-\mu} c^{\mn-\mu} a^{\mu} = (a^{*}a+c^{*}c)^{\mn} = 1,
\eeq
from which both $\qpp$ and $\qpn$ are
idempotents: $(\qpp)^2=\qpp$ and $(\qpn)^2=\qpn$.
\begin{rema}
The coefficients $ \alpha_{n,\mu}$ and $ \beta_{n,\mu}$ above are q-binomial coefficients and their expression is so as to get the identities \eqref{id1} and \eqref{id2}. When computing the q-winding number of the bundles in Sect.~\ref{se:wn} below we shall need  the relation 
\beq\label{remref}
q^{-2\mu + 2\mu (n-\mu)} \beta_{n,\mu} = \alpha_{n,n-\mu} , \qquad \mu = 0, \dots, n \,.
\eeq
which is obtained by a straightforward computation.
\end{rema}

The projections  \eqref{qpro} and \eqref{qpron} play a central role throughout our paper. As a first application of their properties one shows that $\left(\ASU,\Asq,\ca(\U(1))\right)$ is a quantum principal bundle: the relevant proof is in 
App.~\ref{se:pbc}.  

Next, we identify the spaces of equivariant maps $\cl_n$ with the left
$\Asq$-modules of sections $(\Asq)^{n+1}\qpp$ or $(\Asq)^{\mid
n\mid+1}\qpn$, according to whether $n$ is positive or negative. For
this we write any element in the free module $(\Asq)^{\mid n\mid
+1}$ as $\bra{f}=(f_0, f_1,\dots, f_n)$ with $f_\mu \in \Asq$. This
allows one to write equivariant maps as
\begin{align*}
& \phi_f :=\hs{f}{ \Psi^{(n)} }
= \sum\nolimits_{\mu=0}^{n} f_{\mu} \sqrt{\beta_{n,\mu}}\,c^{*\mu} a^{*n-\mu},
\qquad \qquad \mathrm{for} \quad n \geq 0 ,  \\
& \check{\phi}_f :=\hs{f}{ \check{\Psi}^{(n)}  }
= \sum\nolimits_{\mu=0}^{\mn} f_{\mu} \sqrt{\alpha_{n,\mu}}\,c^{\mn-\mu} a^{\mu},
\qquad \qquad \mathrm{for} \quad n \leq 0 . 
\end{align*}
Writing equivariant maps in the above form, the following proposition is easy to establish.
\begin{prop}\label{isoeqsec}
Let $\ce_n:=(\Asq)^{n+1}\qpp$ or $\check{\ce}_n=(\Asq)^{\mn+1}\qpn$ according to whether $n\geq0$ or $n\leq0$.
There are left $\Asq$-modules isomorphisms:
$$
\cl_n ~\xrightarrow{~\simeq~}~ \ce_n, \quad \phi_f  \to \sigma_f := \phi_f \bra{\Psi^{(n)} }
= \bra{f} \qpp ,
$$ with inverse
$$
\ce_n ~\xrightarrow{~\simeq~}~ \cl_n, \quad \sigma_f = \bra{f} \qpp \to \phi_f :=\hs{f}{ \Psi^{(n)} },
$$
and similar maps for the case $n\leq0$.
\end{prop}

%\begin{rema}
%The modules $\cl_n$ could also be realized as right modules via projections. The latter are again the projections $\qpp$ or $\qpn$ above with their role exchanged, that is, $\cl_n\simeq \qpn(\Asq)^{n+1}$ for $n\geq 0$ and  $\cl_n\simeq \qpp(\Asq)^{\mn+1}$ for $n\leq 0$, respectively.
%\end{rema}

\begin{rema}
By exchanging the role of the previous projections $\qpp$ and $\qpn$ one has also an identification of right modules: $\cl_n\simeq \qpn(\Asq)^{n+1}$ and  $\cl_n\simeq \qpp(\Asq)^{\mn+1}$ for $n\geq 0$ and $n\leq 0$ respectively.
\end{rema}

\noindent 
In the sequel, to lighten the notation we shall simply write $\cl_n \simeq
\ce_n:=(\Asq)^{\mn+1}\qpp$ with the positivity or negativity of the label $n$ deciding which projection to be considered, the ones given in \eqref{qpro} and \eqref{qpron} for $n\geq0$ or $n\leq0$ respectively.

%
%\begin{align*}
%&\cl_n \simeq \qpp(\Asq)^{n+1} , \quad \phi_f  \leftrightarrow \sigma_f =
%\bra{f} \qpp = \hs{f}{ \Psi^{(n)} } \bra{\Psi^{(n)} } \qquad \mathrm{for} \quad n > 0 , \nn \\
%&\cl_n \simeq \qpn(\Asq)^{\mn+1} ,  \quad \check{\phi}_f  \leftrightarrow \check{\sigma}_f = \bra{f} \qpn =  \hs{f}{ \check{\Psi}^{(n)} } \bra{\check{\Psi}^{(n)} } \qquad \mathrm{for} \quad n < 0 .
%\end{align*}

\bigskip
The final result that we need in this section is the known
decomposition of the modules $\cl_n$ into representation
spaces under the (right) action of $\su$. 
We have already mentioned the vector space decomposition
$\ASU=\oplus_{n\in\IZ} \cl_n$. 
{}From the definition of the bundles $\cl_n$ in \eqref{libu} and the relations \eqref{relsu} of $\su$ one gets that,
\beq\label{rellb}
E \lt \cl_n \subset \cl_{n+2}, \qquad F \lt \cl_n \subset \cl_{n-2} .
\eeq
On the other hand, commutativity of the left and right actions of $\su$ yields that 
$$
\cl_n \rt  h \subset \cl_n, \qquad \forall \, h\in \su .
$$
In fact, it was already shown in \cite[Thm.~4.1]{SWPod}
that there is also a decomposition, 
\beq
\label{decoln}
\cl_n:=\bigoplus_{J=\tfrac{|n|}{2}, \tfrac{|n|}{2} +1,
\tfrac{|n|}{2} +2, \cdots}V_{J}^{\left(n\right)} ,
\eeq 
with $V_{J}^{\left(n\right)}$ the spin
$J$-re\-pre\-sen\-ta\-tion space (for the right action) of $\su$.
Altogether we get a Peter-Weyl decomposition for $\ASU$ 
(already given in \cite{wor87}).

More explicitly, the
highest weight vector for each $V_{J}^{\left(n\right)}$ in \eqref{decoln} is
$ c^{J-n/2} a^{*J+n/2}  $:
\begin{align}\label{hwv}
& K \lt (c^{J-n/2} a^{*J+n/2}) = q^{n/2} (c^{J-n/2} a^{*J+n/2}), \nn \\
& (c^{J-n/2} a^{*J+n/2}) \rt  K = q^{J} (c^{J-n/2} a^{*J+n/2}) , \qquad 
(c^{J-n/2} a^{*J+n/2}) \rt  F = 0.
\end{align}
Analogously, the lowest weight vector for each $V_{J}^{\left(n\right)}$ in \eqref{decoln} is $a^{J-n/2} c^{*J+n/2}$:
\begin{align*}
& K \lt (a^{J-n/2} c^{*J+n/2}) = q^{n/2} (a^{J-n/2} c^{*J+n/2}),  \\
& (a^{J-n/2} c^{*J+n/2}) \rt  K = q^{-J} (a^{J-n/2} c^{*J+n/2}) , \qquad 
(a^{J-n/2} c^{*J+n/2}) \rt E = 0.
\end{align*}
The elements of the vector spaces $V_{J}^{\left(n\right)}$ can be obtained by acting on the highest weight vectors with the lowering operator 
$\rt  E$, since clearly $\left(c^{J-n/2} a^{*J+n/2}\right)\rt  E\in\cl_{n}$, or explicitly,
$$
K\lt\left[\left(c^{J-n/2} a^{*J+n/2}\right)\rt  E\right]=q^{n/2}
\left[\left(c^{J-n/2} a^{*J+n/2}\right)\rt  E\right] .
$$
To be definite, let us consider $n\geq0$. The first admissible $J$ is $J=n/2$; the highest weight element is $a^{*n}$ and the vector space 
$V_{n/2}^{\left(n\right)}$ is spanned by 
$\{a^{*n}\rt  E^{l}\}$ with $l=0,\dots,n+1$:
$V_{n/2}^{\left(n\right)} = \mathrm{span} \{ a^{*n},c^{*}a^{*n-1},\ldots , c^{*n} \}$. 
Keeping $n$ fixed, the other admissible values of $J$ are $J=s + n/2$ with $s\in\IN$. The vector spaces $V_{s+n/2}^{\left(n\right)}$ are spanned by $\{c^{s}a^{*s+n}\rt  E^{l}\}$ with $l=0,\ldots,2s+n+1$. 
Analogous considerations are valid when $n\leq0$. In this cases, the admissible values of $J$ are $J=s + \mn/2=s-n/2$, the highest weight vector in $V_{s-n/2}^{\left(n\right)}$ is the element $c^{s-n}a^{*s}$, and a basis is given by the action of the lowering operator $\rt  E$, that is  
$V_{s-n/2}^{\left(n\right)}= \mathrm{span} \{ \left(c^{s-n}a^{*s}\right)\rt  E^{l}, 
\; l=0,\ldots,2s-n+1\}$.

%This analysis leads to write an explicit basis for $\cl_{n}$, with $n\in\IZ$ as ($l\in\IN$):
%\beq
%\cl_n=\bigoplus_{s\in\IN}\left(c^{s+(\mn-n)/2}a^{*s+(\mn+n)/2}\right)\rt  E^{l}
%\label{lnba}
%\eeq

We know from (\ref{rellb}) that the left action $F \lt$ maps $\cl_{n}$ 
to $\cl_{n-2}$. If $p\geq0$, the element $a^{*p}$ is the highest weight vector in $V_{p/2}^{(p)}$ and one has that $F\lt a^{*p} \propto c a^{*p-1}$. The element $c a^{*p-1}$ is the highest weight vector in 
$V_{p/2}^{(p-2)}$ since one finds that $(c a^{*p-1}) \rt F=0$ and 
$(c a^{*p-1}) \rt K=q^{p/2} (c a^{*p-1})$. In the same vein, the elements  
$F^{t}\lt a^{*p} \propto c^{t}a^{*p-t}$ are the highest   
weight elements in $V_{p/2}^{(p-2t)} \subset \cl_{p-2t},\,\,t=0,\ldots,p$. Once again, a complete basis of each subspace $V_{p/2}^{(p-2t)}$ is obtained by the right action of the lowering operator $\rt  E$.

With these considerations, the algebra $\ASU$ can be partitioned into finite dimensional blocks which are the analogues of the Wigner D-functions  \cite{Mos} for the group $SU(2)$. To illustrate the meaning of this partition, let us start with the element $a^{*}$, the highest weight vector of the space 
$V_{1/2}^{(1)}$. Representing the left action of $F \lt$ with a horizontal arrow and the right action of  $\rt E$ with a vertical one, yields the box
$$
\begin{array}{ccc}
a^{*} & \to & c \\
\downarrow & & \downarrow \\
-qc^{*} & \to & a 
\end{array}  \;\; ,
$$
where the first column is a basis of the subspace 
$V_{1/2}^{(1)}$, while the second column is a basis of the subspace 
$V_{1/2}^{(-1)}$. Starting from $a^{*2}$ -- the highest weight vector of $V_{1}^{(2)}$ -- one gets:
$$
\begin{array}{ccccc}
a^{*2} & \to & q^{-1/2}\left[2\right]c a^{*} & \to & \left[2\right]c^{2} \\
\downarrow & & \downarrow & & \downarrow \\
-q^{1/2}\left[2\right]c^{*}a^{*} & \to & \left[2\right]\left(a a^{*} - c c^{*}\right) & \to & \left[2\right]^2 q^{1/2}ca \\
\downarrow & & \downarrow & & \downarrow \\
q^{2}\left[2\right]c^{*2} & \to & -q^{3/2}\left[2\right]^2 a c^{*} & \to & q^{3}\left[2\right]^2 a^{2} 
\end{array} \;\; .
$$
The three columns of this box are bases for the subspaces $V_{1}^{(2)}$, $V_{1}^{(0)}$,$V_{1}^{(-2)}$, respectively. The recursive structure is clear. For a positive integer $p$, one has a box $W_{p}$ made up of $\left(p+1\right)\times\left(p+1\right)$ elements. Without explicitly computing the coefficients, one gets:
$$
\begin{array}{ccccccccccc}
a^{*p} & \to & ca^{*p-1} & \to & \ldots & \to & c^{t}a^{*p-t} & \to & \ldots & \to & c^{p} \\
\downarrow & & \downarrow & & \ldots & & \downarrow & & \ldots & & \downarrow \\
c^{*}a^{*p-1} & \to & \ldots & \to & \ldots & \to & \ldots & \to & \ldots & \to & ac^{p-1} \\
\downarrow & & \downarrow & & \ldots & & \downarrow & & \ldots & & \downarrow \\
\ldots & \to & \ldots & \to & \ldots & \to & \ldots & \to & \ldots & \to & \ldots \\
\downarrow & & \downarrow & & \ldots & & \downarrow & & \ldots & & \downarrow \\
c^{*s}a^{*p-s} & \to & \ldots & \to & \ldots & \to & \ldots & \to & \ldots & \to & a^{s}c^{p-s} \\
\downarrow & & \downarrow & & \ldots & & \downarrow & & \ldots & & \downarrow \\
\ldots & \to & \ldots & \to & \ldots & \to & \ldots & \to & \ldots & \to & \ldots \\
\downarrow & & \downarrow & & \ldots & & \downarrow & & \ldots & & \downarrow \\
c^{*p} & \to & ac^{*p-1} & \to & \ldots & \to & a^{t}c^{*p-t} & \to & \ldots & \to & a^{p} \\
\end{array} \;\; .
$$
The space $W_{p}$ is the direct sum of representation spaces for the right action of $\su$, 
$$
W_{p}=\oplus_{t=0}^{p}V_{p/2}^{(p-2t)}, 
$$
and on each $W_{p}$ the quantum Casimir $C_{q}$ act is the same manner from both the right and the left, with eigenvalue (\ref{quca}), that is 
$C_{q}\lt w_{p} = w_{p}\rt  C_{q} = \left([p+\half]^{2}-\tfrac{1}{4}\right)w_{p}$, for all $w_{p}\in W_{p}$.
The Peter-Weyl decomposition for the algebra $\ASU$ is given as,
$$
\ASU = \oplus_{p\in\IN}W_{p} = 
\oplus_{p\in\IN}\left(\oplus_{t=0}^{p}V_{p/2}^{(p-2t)}\right).
$$
At the classical value of the deformation parameter, $q=1$, one recovers the classical Peter-Weyl decomposition for the group algebra of $SU\left(2\right)$.
A straightforward computation shows that the elements in $W_{p}$ which are written (up to a multiplicative constant) as $w_{p:t,r}:=F^{t}\lt a^{*p}\rt  E^{r}$, for $t,r=0,1,\dots,p$, become, in the classical limit, proportional to the Wigner D-functions $D^{p/2}_{t-p/2;r-p/2}$ for $SU(2)$, -- the building blocks \cite{Mos} of the Peter-Weyl decomposition of the algebra of functions on the manifold of the group $SU(2)$.

In order to get elements in the Podle\'s sphere subalgebra $\Asq\simeq\cl_{0}$ out of a highest weight vector $a^{*p}$ we need $p=2l$ to be even and left action of $F^{l}$: $F^{l}\lt a^{*2l} \propto c^{l}a^{*l} \in \Asq$. Then, the right action of $E$ yields a spherical harmonic decomposition,
\beq
\Asq = \oplus_{l \in\IN}V_{l}^{(0)},
\label{decsp}
\eeq
with a basis of $V_{l}^{(0)}$ given by the vectors $F^{l}\lt a^{*2l}\rt  E^{r}$, for $r=0,1,\dots,2l$. At the classical value of the deformation parameter these vectors become the standard spherical harmonics,
$\{Y_{lr}, l\in\IN, r=0,1,\dots,2l \}$, which build up the spherical harmonics decomposition of the algebra of functions on the classical sphere $\mathrm{S}^2$.

\section{The calculi on the quantum principal bundle}\label{se:cotqpb}

The principal bundle $(\ASU,\Asq,\ca(\U(1)))$ is endowed \cite{BM93,BM97} with compatible nonuniversal calculi obtained from the 3-dimensional left-covariant calculus \cite{wor87} on $\SU$
we present first. We then describe the unique left covariant 2-dimensional calculus \cite{pod89} on the sphere $\sq$ obtained by restriction and the projected calculus on the `structure Hopf algebra' $\ca(\U(1))$. All these calculi are compatible in a natural sense.

\subsection{The left-covariant calculus on $\SU$}\label{se:lcc}

The first differential calculus we take on the quantum group $\SU$
is the left-covariant one already developed in \cite{wor87}. It is
three dimensional with corresponding ideal $\mathcal{Q}_{\SU}$
generated by the 6 elements $
\{a^*+q^{2}a-(1+q^{2});c^{2};c^*c;c^{*2};(a-1)c;(a-1)c^*\} $.
%: this calculus turns out to be three dimensional, as
%$dim((Kern\,\varepsilon_{\SU})/\mathcal{Q}_{\SU})=3$.
The quantum tangent space $\mathcal{X}(\SU)$  is generated by the
three elements: 
$$
X_{z} =\frac{1-K^{4}}{1-q^{-2}}, \qquad X_{-}
=q^{-1/2}FK , \qquad X_{+} =q^{1/2}EK , 
$$
whose coproducts are easily found:
$$
\cop X_z = 1\otimes X_z + X_z \otimes K^4, \qquad \cop X_\pm = 1\otimes X_\pm + X_\pm \otimes K^2 .
$$
The dual space of 1-forms $\Omega^1(\SU)$ has a basis \beq
\omega_{z} =a^*\dd a+c^*\dd c , \qquad \omega_{-} =c^*\dd
a^*-qa^*\dd c^*, \qquad \omega_{+} =a\dd c-qc \dd a , \label{q3dom}
\eeq of left-covariant forms, that is
$\Delta_{L}^{(1)}(\omega_{s})=1\otimes\omega_{s}$, with $\Delta_{L}^{(1)}=\Delta^{(1)}$ the (left)
coaction of $\ASU$ onto itself extended to forms. The differential $\dd
: \ASU \to \Omega^1(\SU)$ is 
\beq\label{exts3} 
\dd f =
(X_{+}\triangleright f) \,\omega_{+} + (X_{-}\triangleright f)
\,\omega_{-} + (X_{z}\triangleright f) \,\omega_{z}, 
\eeq 
for all $f\in\ASU$. The above relations \eqref{q3dom} can be inverted to
\begin{align*}
\dd a&=-qc^*\omega_{+}+a\omega_{z} , \qquad
\dd a^*=-q^{2}a^*\omega_{z}+c\omega_{-}\nn\\
\dd c&=a^*\omega_{+}+c\omega_{z} , \qquad
\dd c^*=-q^{2}c^*\omega_{z}-q^{-1}a\omega_{-} ,
\end{align*}
from which one also gets that $\omega_{-}^*=-\omega_{+}$ and
$\omega_{z}^*=-\omega_{z}$. The bimodule structure is:
\begin{align}\label{bi1}
\omega_{z}a=q^{-2}a\omega_{z}, \qquad
\omega_{z}a^*=q^{2}a^*\omega_{z}, \qquad 
\omega_{\pm}a=q^{-1}a\omega_{\pm},\qquad  \omega_{\pm}a^*=qa^*\omega_{\pm} \nn \\
\omega_{z}c=q^{-2}c\omega_{z}, \qquad
\omega_{z}c^*=q^{2}c^*\omega_{z}, \qquad 
\omega_{\pm}c=q^{-1}c\omega_{\pm} , \qquad
\omega_{\pm}c^*=qc^*\omega_{\pm},
\end{align}
Higher dimensional forms can be defined
in a natural way by requiring compatibility for commutation
relations and that $\dd^2=0$. One has:
\beq \label{dformc3} 
\dd
\omega_{z} =-\omega_{-}\wedge\omega_{+} , \quad \dd \omega_{+}
=q^{2}(1+q^{2})\omega_{z}\wedge\omega_{+} , \quad \dd \omega_{-}
=-(1+q^{-2})\omega_{z}\wedge\omega_{-}, 
\eeq 
together with commutation relations:
\begin{align}\label{commc3}
&\omega_{+}\wedge\omega_{+}=\omega_{-}\wedge\omega_{-}=\omega_{z}
\wedge\omega_{z}=0 \nn\\
&\omega_{-}\wedge\omega_{+}+q^{-2}\omega_{+}\wedge\omega_{-}=0
\nn\\
& \omega_{z}\wedge\omega_{-}+q^{4}\omega_{-}\wedge\omega_{z}=0, \qquad
\omega_{z}\wedge\omega_{+}+q^{-4}\omega_{+}\wedge\omega_{z}=0.  \end{align}
Finally, there is a unique top form $\omega_{-}\wedge\omega_{+}\wedge\omega_{z}$.

\subsection{The standard calculus on $\sq$}\label{se:cals2}

The restriction of the above 3D calculus to the sphere $\sq$ yields
the unique left covariant 2-dimensional calculus on the latter \cite{maj05}. An
evolution of this approach has led \cite{SW04} to a description of
the unique 2D calculus of $\sq$ in term of a Dirac operator. The `cotangent bundle' $\Omega^1(\sq)$ is shown to be isomorphic
to the direct sum $\cl_{-2}\oplus\cl_2$, that is the line bundles
with winding number $\pm 2$. Since the element $K$ acts as the
identity on $\Asq$, the differential \eqref{exts3} becomes, when restricted to the latter, 
\begin{align*}
\dd f &= (X_{-}\triangleright f) \,\omega_{-} + (X_{+}\triangleright f) \,\omega_{+}  \\
& =  (F \triangleright f) \,
(q^{-1/2}\omega_{-}) + (E \triangleright f) \,(q^{1/2}\omega_{+}) ,  \qquad \mathrm{for} \, f\in\Asq .
\end{align*}
These leads to break the exterior differential into a holomorphic and an
anti-holomorphic part, $\dd = \delb + \del$, with:
\begin{align*}
&\delb f=\left(X_{-}\triangleright f\right)\omega_{-} = (F \triangleright f) \,
(q^{-1/2}\omega_{-}) , \nn \\
&\del f=\left(X_{+}\triangleright f\right)\omega_{+} = (E \triangleright f)
\,(q^{1/2}\omega_{+}) , \qquad \mathrm{for} \quad f\in\Asq .
\end{align*}
An explicit computation on the generators \eqref{podgens} of $\sq$
yields: 
%\beq 
%\delb\bm=\frac{1}{\left(1+q^{2}\right)^{1/2}}\,a^{2} \,
%\omega_{-}, \qquad \delb\bz=q \, ca  \, \omega_{-}, \qquad
%\delb\bp=\frac{q^{2}}{\left(1+q^{2}\right)^{1/2}}\,c^{2} \,
%\omega_{-} . \label{dellbb} 
%\eeq and:
% \beq
%\del\bp=\frac{q^{3}}{\left(1+q^{2}\right)^{1/2}}\,a^{*2} \,
%\omega_{+}, \qquad \del\bz=-q^{2}\,c^{*}a^{*} \, \omega_{+}, \qquad
%\del\bm=\frac{q^{3}}{\left(1+q^{2}\right)^{1/2}}\,c^{*2} \,
%\omega_{+} \label{dellb} 
%. \eeq 
$$
\begin{array}{lll}
\delb \bm= q^{-1}\, a^{2} \, \omega_{-}, 
& \delb\bz=q \, ca  \, \omega_{-}, 
& \delb\bp= q\, c^{2} \, \omega_{-} , \\
~\\
\del\bp= q^{2} \,a^{*2} \, \omega_{+}, 
& \del\bz= -q^{2}\, c^{*}a^{*} \, \omega_{+}, 
& \del\bm= q^{2}\, c^{*2} \, \omega_{+} .
\end{array}
$$
The above shows that:
$\Omega^1({\sq})=\Omega^{1}_{-}(\sq) \oplus\Omega^{1}_{+}(\sq)$ where
$\Omega^{1}_{-}(\sq)\simeq \cl_{-2} \simeq \delb(\Asq)$ is the
$\Asq$-bimodule generated by: 
$$
\{\delb\bm,\delb\bz,\delb\bp\}=\{a^{2},ca,c^{2}\}\,\omega_{-} =
q^{2}\omega_{-}\{a^{2},ca,c^{2}\} 
$$
and $\Omega^{1}_{+}(\sq)\simeq
\cl_{+2}\simeq \del(\Asq)$ is the one generated by:
$$
\{\del\bp,\del\bz,\del\bm\}=\{a^{*2},c^{*}a^{*},c^{*2}\}\,
\omega_{+}= q^{-2} \omega_{+} \{a^{*2},c^{*}a^{*},c^{*2}\} .
$$
That these two modules of forms are not free is also expressed by
 the existence of relations among the differential:
%\begin{align*}
%&\del\bz=\left(q^{-1}+q^{-3}\right)\bm\del\bp-\left(q+q^{3}\right)\bp\del\bm , \\
%&\delb\bz=\left(q^{-1}+q\right)\bp\delb\bm-\left(q^{-5}+q^{-3}\right)\bm\delb\bp .
%\end{align*}
$$
\del\bz= q^{-1} \bm\del\bp - q^{3} \bp\del\bm , \qquad 
\delb\bz= q \bp\delb\bm - q^{-3} \bm\delb\bp .
$$
The 2D calculus on $\sq$ has a unique top 2-form $\omega$ with
$\omega f = f \omega$,  for all $f\in\Asq$ and $\Omega^2({\sq})$ is
the free $\Asq$-module generated by $\omega$, that is
$\Omega^2({\sq})=\omega \Asq = \Asq \omega$. Now, both
$\omega_{\pm}$ commutes with elements of $\Asq$ and so does
$\omega_{-}\wedge\omega_{+}$, which is taken as the natural
generator $\omega=\omega_{-}\wedge\omega_{+}$ of $\Omega^2({\sq})$. Writing  any 1-form as $\alpha = x
\omega_{-} + y \omega_{+} \in \cl_{-2} \omega_{-} \oplus \cl_{+2}
\omega_{+}$, the product of 1-forms is:
$$
(x \omega_{-} + y \omega_{+} ) \wedge (t \omega_{-} + z \omega_{+} ) =
(q^{-2} y t - x z)  \omega_{+}\wedge\omega_{-}.
$$
{}From \eqref{dformc3} it is natural to ask that $\dd \omega_-=\dd
\omega_+=0$ when restricted to $\sq$. Then, the exterior derivative
of any 1-form  $\alpha = x \omega_{-} + y \omega_{+} \in \cl_{-2}
\omega_{-} \oplus \cl_{+2} \omega_{+}$ is  given by:
\begin{align}\label{d1f}
\dd \alpha & = \dd (x \omega_{-} + y \omega_{+}) = \del x \wedge \omega_{-} + \delb y \wedge \omega_{+} \nn \\
&= (X_+ \lt x  - q^{-2}  X_- \lt y) \, \omega_{+}\wedge\omega_{-} =
q^{-1/2} ( E \lt x  - q^{-1} F \lt y ) \, \omega_{+}\wedge\omega_{-} ,
\end{align}
since $K\lt$ acts as $q^{\mp}$ on $\cl_{\mp 2}$. Notice that in the
above equality, both $E \lt x$ and $F \lt y$ belong to $\Asq$, as it
should be. We summarize the above results in the following
proposition.
\begin{prop}\label{2dsph}
The 2D differential calculus on the sphere $\sq$ is given by:
$$
\Omega^{\bullet}({\sq}) = \Asq \oplus \left(\cl_{-2} \omega_{-}
\oplus \cl_{+2} \omega_{+} \right) \oplus \Asq \omega_{+}\wedge\omega_{-} ,
$$
with multiplication rule
$$
\big(f_0; x,y;f_2\big)\big(g_0;t,z;g_2\big)=\big(f_0g_0; f_0t+xg_0,f_0z+yg_0;
f_0g_2+f_2g_0+q^{-2}y t - x z\big),
$$
and  exterior differential $\dd = \delb + \del$:
$$
\begin{array}{ll}
f \mapsto (q^{-1/2} F\lt f, q^{1/2} E \lt f ), & \mathrm{for} \, f\in\Asq , \\
~ & \\
(x,y) \mapsto q^{-1/2} ( E \lt x  - q^{-1} F \lt y ) , & \mathrm{for} \,
(x,y)\in\cl_{-2} \oplus \cl_{+2} .
\end{array}
$$
\end{prop}
\begin{rema}
It is evident that  the relevant operators for the calculus on
$\Asq$ are $\{E,F\}$ rather than $\{X_+,X_-\}$; indeed, it is proven
in \cite{SW04} that $\{E,F\}$ span the quantum tangent space of the
2D calculus $(\Omega^2({\sq}),\dd)$. We mantain $\{X_+,X_-\}$ as
well since these are the operators that will be lifted to the total
space $\ASU$ via the connection and that will enter the gauged
Laplacian later on.
\end{rema}

\subsection{The calculus on the structure group}\label{se:csg}

{}From App.~\ref{se:qpb} we know that a quantum principal bundle
with nonuniversal calculi comes with compatible calculi on the `total space
algebra' $\mathcal{P}$ and on the `structure Hopf algebra'
$\ch$. The compatibility determines a calculus on
the `base space algebra' $\mathcal{B}$.

For the Hopf bundle over the sphere $\sq$ we have already given the
calculus on $\sq$ starting from the left covariant one on $\SU$.
Here we give the calculus on $\U(1)$ which makes them compatible from
the quantum principal bundle point of view.

The strategy \cite{BM93} consists in defining the calculus on $\U(1)$
via the Hopf projection $\pi$ in \eqref{qprp}. In particular, out of
the $\mathcal{Q}_{\SU}$ which determines the left covariant calculus
on $\SU$, one defines a right ideal
$\mathcal{Q}_{\U(1)}=\pi(\mathcal{Q}_{\SU})$ for the calculus on
$\U(1)$. One finds that $\mathcal{Q}_{\U(1)}$ is generated by the
element $\{z^{*}+q^{2}z-(1+q^{2})\}$ and the corresponding calculus
is 1-dimensional, and bicovariant. %%
%\footnote{This differential calculus is slightly
%different from that defined in \cite{brz94} because their
%projection is defined not on $\SU$, but on
%$SO_{q}(3)=SU_{q^{2}}(2)/\zed^{2}$.}.
%%
Its quantum tangent space is generated by  
\beq\label{vvf} 
X=X_{z} =\frac{1-K^{4}}{1-q^{-2}},
\eeq 
with dual 1-form given by $\omega_{z}$. Explicitly, one finds that
$$ 
\omega_{z} = z^{*} \dd z , \qquad 
\dd z=z\omega_{z} , \qquad
\dd z^{*} =-q^{2}z^{*} \omega_{z} , 
$$
and noncommutative commutation relations,
$$
\omega_{z} z = q^{-2} z\omega_{z}, \qquad
\omega_{z} z^{*} = q^{2}z^{*}\omega_{z} , \qquad
z \dd z=q^{2}(\dd z) z .
$$

We know from Sect.~\ref{se:avb} and App.~\ref{se:pbc} that the datum
$(\ASU,\Asq,\ca(\U(1)))$ is a `topological' quantum principal bundle. We have now a specific differential calculus both on the total space $\ASU$ (the 3D left covariant calculus) and on $\ca(\U(1))$ (obtained from it via the same projection $\pi$ in \eqref{qprp} giving the bundle structure). Moreover, from the calculus on  $\ASU$, we also obtained via restriction a calculus on the base space $\Asq$. 
It remains to show that these calculi are compatible so that the datum $(\ASU,\Asq,\ca(\U(1)); \cq_{\SU}, \cq_{\ca(\U(1))})$ is a quantum principal bundle with nonuniversal calculi. As a result, the vector field $X_z$ is vertical. 
The details showing the compatibility are in App.~\ref{se:cc}. In particular, from the analysis of the appendix, the vector field \eqref{vvf} is a `vertical' vector field for the fibration.

\section{The monopole connection and its curvature}\label{se:con}

A connection on the quantum principal bundle with respect to the left covariant calculus 
$\Omega(\sq)$ will be the crucial ingredient for the gauged Laplacian operator. 
The connection will determine a covariant derivative on the module of equivariant maps $\cl_n$ which will, in turn, be shown to correspond to the Grassmann connection on the modules $(\Asq)^{\mn+1}\qpp$. We shall also derive corresponding expressions for the curvature of the connection. 

\subsection{Enter the connection}

The most efficient way to define a connection on a quantum principal
bundle (with given calculi) is by splitting the 1-forms on the total
space into horizontal and vertical ones \cite{BM93,BM97}. Since
horizontal 1-forms 
%$\Omega_{\mathrm{hor}}{(\mathcal{P})}$ 
are given
in the structure of the principal bundle, one needs a projection on
forms whose range is the subspace of vertical ones. The projection
is required to be covariant with respect to the right coaction of
the structure Hopf algebra.

For the principal bundle over the quantum sphere $\sq$ that we are
considering, a principal connection is a covariant left module 
projection $\Pi : \Omega^1(\SU) \to \Omega^1_{\mathrm{ver}}(\SU)$, that is
$\Pi^2=\Pi$ and $\Pi(x \alpha) = x \Pi(\alpha)$, for
$\alpha\in\Omega^1(\SU)$ and $x \in\mathcal{A}(\SU)$;  equivalently
it is a covariant splitting $\Omega^1(\SU)=\Omega^1_{\mathrm{ver}}(\SU)
\oplus\Omega^1_{\mathrm{hor}}(\SU)$. The covariance of the connection
is the requirement that 
$$
\Delta_{R}^{(1)}\,\Pi=\left(\Pi\otimes\id\right)\circ \Delta_{R}^{(1)} , 
$$ 
with $\Delta_{R}^{(1)}$ the extension to 1-forms of the coaction $\Delta_{R}$ (\ref{cancoa}) of the structure Hopf algebra. It
is not difficult to realize that with the left covariant 3D calculus on $\ASU$, a
basis for $\Omega^1_{\mathrm{hor}}(\SU)$ is given by $\omega_{-},
\omega_{+}$ (a proof is in the App.~\ref{se:cc}). Furthermore:
$$
\Delta_{R}^{(1)}(\omega_{z})=\omega_{z}\otimes 1, \qquad
\Delta_{R}^{(1)}(\omega_{-})=\omega_{-}\otimes z^{*2}, \qquad
\Delta_{R}^{(1)}(\omega_{+})=\omega_{+}\otimes z^{2},
$$
and a natural choice of a connection is to define $\omega_{z}$ to be
vertical \cite{BM93,maj05}: 
$$
\Pi_{z}(\omega_{z}):=\omega_{z}, \qquad
\Pi_{z}\left(\omega_{\pm}\right):=0. 
$$
With a connection, one has
a covariant derivative of equivariant maps, 
$$
\nabla:\ce \to \Omega(\sq) \otimes_{\Asq} \ce, 
\qquad \nabla := (\id - \Pi_{z}) \circ \dd , 
$$
and one readily shows the Leibniz rule
property: $\nabla(f \phi)=f \nabla(\phi)+(\dd f) \otimes \phi$, for all
$\phi\in\ce$ and $f\in\Asq$. We shall take for $\ce$ the line
bundles $\cl_n$ of \eqref{libu}. Then, with the left covariant 2D
calculus on $\Asq$ (coming from the left covariant 3D calculus on
$\ASU$ as explained in Sect.~\ref{se:cals2}) we have,
\begin{align}\label{coder2d}
\nabla \phi &= \left(X_{+}\triangleright\phi\right)\omega_{+}
+\left(X_{-}\triangleright\phi\right)\omega_{-} \nn \\
&= q^{-n-2} \omega_{+} \left(X_{+}\triangleright\phi\right)
+ q^{-n+2} \omega_{-} \left(X_{-}\triangleright\phi\right),
\end{align}
since $X_{\pm} \lt \phi \in \cl_{n\pm2}$. 
The left action of $\su$ on $\Omega(\sq)$ is defined by requiring that it commutes with the exterior derivative $\dd$ while on $\Omega(\sq) \otimes_{\Asq}\ce$ is defined in a obvious way as $h\lt (\alpha \otimes \phi) = \left(\co{h}{1} \lt \alpha\right) \otimes \left(\co{h}{2}\phi \right)$, for $h\in\su$, with notation $\Delta(h)=\co{h}{1} \otimes {\co{h}{2}}$ for the coproduct. 
Then, while
$X_{\pm} \lt \phi \notin \cl_n$, one checks that 
$K\lt \left(\nabla\phi\right) = q^{n} \nabla\phi$ and $\nabla\phi\in\Omega(\sq) \otimes_{\Asq} \ce$ as it should be.

It was shown in \cite{HM98} that with the universal
calculi on the principal bundle, the module $\ce$ is projective and
the covariant derivative $\nabla$ corresponds to the Grassmann
connection of the corresponding projection. We have a similar result
for the left covariant calculi when taking as modules $\ce$ the
line bundles $\cl_n$ of \eqref{libu}.

\begin{prop}
Let $\ce$ be the line bundle $\cl_n$ defined in \eqref{libu}. With
the left $\Asq$-modules isomorphism $\cl_n \simeq
\ce_n:=(\Asq)^{\mn+1}\qpp$ of Prop.~\ref{isoeqsec} extended in a
natural way to forms $\Omega(\sq) \otimes_{\Asq} \cl_n \simeq \Omega(\sq)
\otimes_{\Asq} \ce_n$, the covariant derivative on $\cl_n$
corresponds to the Grassmann connection on $\ce_n$, that is,
$$
 ( \dd \sigma_{\phi} ) \, \qpp = \sigma_{\nabla\phi} ,
$$
for $\phi\in\cl_n$ and corresponding section $\sigma_{\phi}\in\ce_n$.
\begin{proof}

We consider only the case $n\geq0$ since the proof for the case $n\leq0$
is the same. {}From Prop.\ref{isoeqsec}, to any $\phi\in\cl_n$ there
corresponds the section $\sigma_\phi= \phi \bra{\Psi^{(n)}} \in
\ce_n$, having components $(\sigma_{\phi})_{\nu} = (\phi\,
a^{n-\nu}c^{\nu}) \sqrt{\beta_{n,\nu}}$, $\nu=0,1, \dots, n$.
Similarly, to the covariant derivative $\nabla\phi\in\Omega(\sq)
\otimes_{\Asq} \cl_n$ we associate the $\Omega(\sq)$-valued section
$\sigma_{\nabla\phi}= (\nabla\phi) \bra{\Psi^{(n)}} \in  \Omega(\sq)
\otimes_{\Asq} \ce_n$. Explicitly:
\begin{align*}
\sigma_{\nabla\phi}
& =\left[ \left(X_{+}\triangleright\phi\right) \omega_{+}
+ \left(X_{-}\triangleright\phi \right)\omega_{-} \right] \bra{\Psi^{(n)}}
\nn\\
& =
q^{-n} \left[ \left(X_{+}\triangleright\phi\right) \bra{\Psi^{(n)}}\, \right] \omega_{+}
+ q^{-n} \left[ \left(X_{-}\triangleright\phi \right) \bra{\Psi^{(n)}}\, \right]
\omega_{-}  ,
\end{align*}
since $\bra{\Psi^{(n)}}\in\cl_{-n}$ and $\omega_{\pm}\cl_{-n}
\subset q^{-n}\cl_{-n} \omega_{\pm}$, from the commutation relations
in the second column in \eqref{bi1}. In components: 
$$
\left[\sigma_{\nabla\phi}\right]_\mu = q^{-n} \left(
X_{+}\triangleright \phi \right)\, (a^{n-\mu}c^{\mu})
\sqrt{\beta_{n,\mu}}\, \omega_{+} + q^{-n} \left(X_{-}\triangleright
\phi \right)\, (a^{n-\mu} c^{\mu}) \sqrt{\beta_{n,\mu}}\,
\omega_{-}. 
$$
On the other hand, for the Grassmann connection
$\nabla \sigma_{\phi} := \dd (\sigma_{\phi} )\, \qpp$ acting on the section $\sigma_{\phi}$,
using the commutativity of  $\omega_{\pm}$ with $\cl_0=\Asq$, one
finds:

\begin{align*}
& \left[ \dd \left(\sigma_{\phi}\right) \qpp \right]_{\mu}
= \sum\nolimits_{\nu=0}^{n} \dd \left(\sigma_{\phi}\right)_{\nu} \qpp_{\nu\mu}
= \sum\nolimits_{\nu=0}^{n} \dd \left( \phi\,  a^{n-\nu}c^{\nu} \right)
\sqrt{\beta_{n,\nu}}\, \qpp_{\nu\mu}  \\
& \qquad  = \sum\nolimits_{\nu=0}^{n}
X_{+}\triangleright \left( \phi\,  a^{n-\nu}c^{\nu} \right)
\sqrt{\beta_{n,\nu}}\, \qpp_{\nu\mu}\, \omega_{+} +
\sum\nolimits_{\nu=0}^{n} X_{-}\triangleright \left(\phi\, a^{n-\nu} c^{\nu}\right)
\, \sqrt{\beta_{n,\nu}}\, \qpp_{\nu\mu}\, \omega_{-} \\
& \qquad  = \sum\nolimits_{\nu=0}^{n}
\left[  \phi\,  X_{+}\triangleright \left(a^{n-\nu}c^{\nu} \right)
+ q^{-n} \left( X_{+}\triangleright \phi \right)\, (a^{n-\nu}c^{\nu}) \right] \,
\sqrt{\beta_{n,\nu}}\, \qpp_{\nu\mu}\, \omega_{+} \\
& \qquad\qquad \qquad\qquad\qquad\qquad\qquad\qquad
+ \sum\nolimits_{\nu=0}^{n} q^{-n} \left(X_{-}\triangleright \phi\left)\, \right(a^{n-\nu} c^{\nu}\right)
\, \sqrt{\beta_{n,\nu}}\, \qpp_{\nu\mu}\, \omega_{-} \\
&\qquad  =
q^{-n} \left( X_{+}\triangleright \phi \right)\,
\sum\nolimits_{\nu=0}^{n}
(a^{n-\nu}c^{\nu}) \sqrt{\beta_{n,\nu}}\, \qpp_{\nu\mu}\, \omega_{+}
\\
& \qquad\qquad \qquad\qquad\qquad\qquad\qquad\qquad
+ q^{-n} \left(X_{-}\triangleright \phi \right)\,
\sum\nolimits_{\nu=0}^{n} (a^{n-\nu} c^{\nu}) \sqrt{\beta_{n,\nu}}\, \qpp_{\nu\mu}\, \omega_{-} \\
& \qquad  =
q^{-n} \left( X_{+}\triangleright \phi \right)\,
(a^{n-\mu}c^{\mu}) \sqrt{\beta_{n,\mu}}\, \omega_{+}
+ q^{-n} \left(X_{-}\triangleright \phi \right)\,
(a^{n-\mu} c^{\mu}) \sqrt{\beta_{n,\mu}}\, \omega_{-}  =\left[\sigma_{\nabla\phi}\right]_\mu .
 \end{align*}
Here the fourth equality follows from the vanishing
$$\sum\nolimits_{\nu=0}^{n}
X_{+}\triangleright \left(a^{n-\nu}c^{\nu} \right) 
\sqrt{\beta_{n,\nu}}\, \qpp_{\nu\mu}\, = q^{-n} [ X_{+} \triangleright  (1) ] = 0.
$$
As mentioned, in a similar fashion one proves the same result for the case $n\leq0$.
\end{proof}
\end{prop}

\subsection{The curvature}

Having the connection we can work out an explicit expression for its curvature, the $\Asq$-linear map $\nabla^2: \ce \to \Omega^2(\sq) \otimes_{\Asq} \ce$, by definition.
\begin{prop}
Let $\ce$ be the line bundle $\cl_n$ defined in \eqref{libu} endowed with
the connection $\nabla$, for  the canonical left covariant 2D
calculus on $\Asq$, given in \eqref{coder2d}. 
Then, with $\phi\in\cl_n$, its curvature is given by: 
$$
\nabla^2 \phi = - q^{-2n-2} \omega_+\wedge\omega_-  \left(X_z \lt \phi\right) \, , 
$$
with $X_z$ the vertical vector field in \eqref{vvf}. As an element in 
$\mathrm{Hom}_{\Asq}(\cl_n, \Omega(\sq) \otimes_{\Asq} \cl_n)$ one finds:
$$
\nabla^2 =q^{-n-1} [n] \, \omega_+\wedge\omega_- \, .
$$
\begin{proof}
Using \eqref{coder2d} and being $\dd\omega_-=\dd\omega_+=0$ on $\sq$, we have
\begin{align*}
\nabla( \nabla \phi) &= - q^{-n-2} \omega_{+} \wedge \left(X_{-} X_{+}\triangleright\phi\right) \omega_-
-q^{-n+2} \omega_{-} \wedge \left(X_{+} X_{-}\triangleright\phi\right))\omega_{+} \\
& = - q^{-2n-2} \omega_{+}\wedge\omega_- \left(X_{-} X_{+}\triangleright\phi\right)
-q^{-2n+2} \omega_{-} \wedge \omega_{+} \left(
X_{+} X_{-}\triangleright\phi\right) \\
& = -q^{-2n-2} \omega_{+}\wedge \omega_-  
\left( X_{-} X_{+} - q^{2} X_{+} X_{-} \right) \triangleright\phi ,
\end{align*}
and the first statement follows from the relation 
$X_{-} X_{+} - q^{2} X_{+} X_{-} = X_z$. The second statement comes from the computation of  $\left(X_z \lt \phi\right)$ on $\phi\in\cl_n$. Since $X_z \lt (\Asq) = 0$, it is evident that the curvature is $\Asq$-linear. 
\end{proof}
\end{prop}
On the other hand, from the action of the connection as $ \nabla \sigma = ( \dd \sigma) \, \qpp$ on the module $\ce_n:=(\Asq)^{\mn+1}\qpp$, a straightforward computation gives for its curvature $F_\nabla = \nabla^2$ as an element in $\Omega(\sq) \otimes_{\Asq} \ce_n$ the expression:
\beq\label{curvsec}
F_\nabla = - \dd \qpp \wedge\, \dd \qpp\,  \qpp \, .
\eeq
To compare these two expressions for the curvature, we need an intermediate result.
\begin{lemm}\label{lemproj}
Let $\qpp$ denote the projection given in \eqref{qpro} and \eqref{qpron}, 
for $n\geq0$ or $n\leq0$, respectively. With
the standard 2D calculus on $\sq$ of Sect.~\ref{se:cals2} one finds: 
\begin{align*}
& \dd \qpp \wedge\, \dd \qpp\,  \qpp\, = -q^{-n-1} [n]
~\qpp\, \omega_{+}\wedge\omega_{-} , \\
& \qpp\,\dd \qpp \wedge\, \dd \qpp\, = -q^{-n-1} [n]
 ~\qpp\, \omega_{+}\wedge\omega_{-} 
.
\end{align*}
%$$
%\begin{array}{ll}
% \qpp\,\dd \qpp \wedge\, \dd \qpp\, = \,- q^{n+1} [n]
% ~\qpp\, \omega_{-}\wedge\omega_{+}  &\qquad \mathrm{for} \,n\geq0 , \\
% & ~ \\
% \qpn\,\dd \qpn \wedge\, \dd \qpn\, = \,  q^{n+1} [n]
% ~\qpn\,
%\omega_{-}\wedge\omega_{+} &\qquad \mathrm{for} \,n\geq0 .
%\end{array}
%$$
\begin{proof}
This is proved by explicit computation. We explicitly consider only
the case $n\geq0$ since the proof for the case $n\leq0$ is the same. In
components:
%$$
%\qpp=\sqrt{\beta_{n,\mu}\beta_{n,\nu}}\,
%c^{*\mu}a^{*n-\mu}a^{n-\nu}c^{\nu}
%$$
\begin{multline*}
\left(\qpp \dd \qpp\right)_{\sigma\nu}=
\sum\nolimits_{\mu=0}^{n}\qpp_{\sigma\mu}d\qpp_{\mu\nu} \\=
\sqrt{\beta_{n,\sigma}\beta_{n,\nu}}\,c^{*\sigma}a^{*n-\sigma}
\sum\nolimits_{\mu=0}^{n}\beta_{n,\mu}a^{ n-\mu}c^{
\mu}\, \dd \left[c^{*\mu}a^{*n-\mu}a^{ n-\nu}c^{ \nu}\right].
\end{multline*}
In the above the anti-holomorphic part vanishes:
\begin{align*}
& \left(\qpp \delb \qpp\right)_{\sigma\nu} =
\sqrt{\beta_{n,\sigma}\beta_{n,\nu}}\,c^{*\sigma}a^{*n-\sigma}
\sum\nolimits_{\mu=0}^{n}\beta_{n,\mu}a^{ n-\mu}c^{\mu} \left(F \triangleright\left[c^{*\mu}a^{*n-\mu}a^{n-\nu}
c^{\nu}\right]\right)\,q^{-1/2} \omega_{-} 
\\
& \qquad =\sqrt{\beta_{n,\sigma}\beta_{n,\nu}}\,c^{*\sigma}a^{*n-\sigma}
\sum\nolimits_{\mu=0}^{n}\beta_{n,\mu}a^{ n-\mu}c^{
\mu} \left(F \triangleright \left[
c^{*\mu}a^{*n-\mu}\right] \right) q^{-n/2}a^{ n-\nu}c^{
\nu}\,q^{-1/2} \omega_{-}
\\
& \qquad =\sqrt{\beta_{n,\sigma}\beta_{n,\nu}}\,c^{*\sigma}a^{*n-\sigma}
\left[\sum\nolimits_{\mu=0}^{n}\beta_{n,\mu} F \triangleright\left(
a^{ n-\mu}c^{
\mu}c^{*\mu}a^{*n-\mu}\right)\right]q^{-n}a^{ n-\nu}c^{\nu} \,
q^{-1/2} \omega_{-} 
\\
& \qquad =\sqrt{\beta_{n,\sigma}\beta_{n,\nu}}\,c^{*\sigma}a^{*n-\sigma}
\left[ F \triangleright \left( 1 \right) \right] q^{-n}a^{ n-\nu}c^{\nu} \,
q^{-1/2} \omega_{-}=0,
\end{align*}
where the second and the third equality  follow from the zero action
of $F \lt$ on any power of $a$ or $c$.
As for the holomorphic part:
\begin{align*}
& \left(\qpp \del \qpp\right)_{\sigma\nu} = \sqrt{\beta_{n,\sigma}\beta_{n,\nu}}\,c^{*\sigma}a^{*n-\sigma}
\sum\nolimits_{\mu=0}^{n}\beta_{n,\mu}a^{n-\mu}c^{
\mu}\left[E\triangleright\left(c^{*\mu}a^{*n-\mu}a^{
n-\nu}c^{\nu}\right)\right]
q^{1/2} \omega_{+} 
\\
& \qquad =\sqrt{\beta_{n,\sigma}\beta_{n,\nu}}\,c^{*\sigma}a^{*n-\sigma}
\sum\nolimits_{\mu=0}^{n}\beta_{n,\mu}a^{ n-\mu}c^{
\mu}c^{*\mu}a^{*n-\mu}
\left[E\triangleright\left(a^{n-\nu}c^{\nu}\right)\right]
q^{-n/2} q^{1/2} \omega_{+} \\ 
& \qquad = \sqrt{\beta_{n,\sigma}\beta_{n,\nu}}\,c^{*\sigma}a^{*n-\sigma}
 \left( 1 \right)  \left[E\triangleright\left(a^{n-\nu}c^{
\nu}\right)\right]
q^{-n/2} q^{1/2} \omega_{+}
\\ 
& \qquad =
\left[E\triangleright \sqrt{\beta_{n,\sigma}\beta_{n,\nu}}\,c^{*\sigma}a^{*n-\sigma}
\left(a^{ n-\nu}c^{
\nu}\right)\right]
q^{1/2} \omega_{+} = \left(E\triangleright\qpp_{\sigma\nu}\right)
q^{1/2} \omega_{+} ,
\end{align*}
using now the zero action of $E \lt$ on any power of $a^*$ or $c^*$.
Having established that: 
$$
\qpp \dd \qpp = \qpp \del \qpp =
\left(E \triangleright \qpp \right) q^{1/2} \omega_{+}, 
$$
we have in turn:
\begin{align*}
\qpp\,\dd \qpp\wedge \dd \,\qpp &= \qpp\, \del \qpp\wedge \delb \,\qpp \\ & =
\left(E\triangleright\qpp\right) \omega_{+}\wedge \left(F\triangleright\qpp\right)
\omega_{-} = \left(E\triangleright\qpp\right)\left(F\triangleright\qpp\right)
q^{-2}\omega_{+}\wedge\omega_{-} .
\end{align*}
The last equality is easily found: being $\qpp_{\mu\sigma}\in\cl_0$
one has $F\triangleright \qpp_{\mu\sigma}\in\cl_{-2}$ (see
\eqref{rellb}), but $\omega_{+}\cl_{-2} \subset q^{-2}\cl_{-2}
\omega_{+}$,  from the second column in the commutation relations
\eqref{bi1}. 
We need to compute:
\begin{multline*}
\left[\left(E \triangleright\qpp\right)\left(F\triangleright\qpp\right)\right]_{\mu\nu} = \\
= q^{-n}\, \sqrt{\beta_{n,\mu}\beta_{n,\nu}}\,c^{*\mu}a^{*n-\mu}
\left\{\sum\nolimits_{\sigma=0}^{n}\beta_{n,\sigma}\left[E\triangleright\left(a^{
n-\sigma}c^{\sigma}\right)\right]\left[F\triangleright\left(c^{*\sigma}a^{*n-\sigma}\right)\right]\right\}
a^{n-\nu}c^{\nu} .
\end{multline*}
For the term in the curly brackets we use the twisted derivation
property of $E$ in:
\begin{multline*}
 0=E \triangleright\left[
\left(a^{n-\sigma}c^{\sigma}\right) F \triangleright\left(c^{*\sigma}a^{*n-\sigma} \right) \right]
= \\ =
q^{(n-2)/2}\left[E\triangleright\left(a^{n-\sigma}c^{\sigma}\right)\right]
\left[F\triangleright\left(c^{*\sigma}a^{*n-\sigma}\right)\right] +
q^{n/2} \left(a^{n-\sigma}c^{\sigma}\right) \left[E F \triangleright
c^{*\sigma}a^{*n-\sigma}\right] .
\end{multline*}
Then, rearranging the terms and summing over $\sigma$ one arrives at
\begin{multline*}
\left[ \left(E\triangleright\qpp\right)\left(F\triangleright\qpp\right) \right]_{\mu\nu} = \\ =
-q^{-n+1} \sqrt{\beta_{n,\mu}\beta_{n,\nu}}\,c^{*\mu}a^{*n-\mu}
\left\{\sum\nolimits_{\sigma=0}^{n}\beta_{n,\sigma}a^{
n-\sigma}c^{\sigma} \left[ E F \triangleright\left( c^{*\sigma}a^{*n-\sigma} \right) \right] \right\}
a^{n-\nu}c^{\nu} .
\end{multline*}
With a direct computation one gets that: 
$E F \triangleright\left(c^{*\sigma} a^{*n-\sigma}\right) = 
[n]\, c^{*\sigma} a^{*n-\sigma}$,  
giving in turn, $\left[ \left(E\triangleright\qpp\right)\left(F\triangleright\qpp\right) \right]_{\mu\nu} = -q^{-n+1} [n] \,
\qpp_{\mu\nu}$, 
and one finally obtains the claimed result:
$$
\left(\qpp\,\dd \qpp\wedge \dd\,\qpp\right)_{\mu\nu}=-q^{-n-1} [n] \,
\qpp_{\mu\nu}\, \omega_{+}\wedge\omega_{-}.
$$
The case $n\leq0$ goes in a similar fashion. One has: 
$$
\qpp\,\dd
\qpp\,=\qpp\,\delb \qpp\,=\left(F\triangleright\qpp\right) q^{-1/2}
\omega_{-} , 
$$
leading to
\begin{align*}
\qpp\,\dd \qpp\wedge\dd \qpp &=\left(F\triangleright\qpp\right)\omega_{-}\wedge
\left(E\triangleright\qpp\right)\omega_{+}=
-q^{\mn+1} [\mn]\, \qpp\, \omega_{-}\wedge\omega_{+} \\
& = -q^{-n-1} [n]\, \qpp\, \omega_{+}\wedge\omega_{-}.
\end{align*}
The second equality in the Lemma is proved in a similar fashion,
using 
$$ 
\left(\dd \qpp\right) \qpp = \left(\delb
\qpp \right) \qpp= \left(F \triangleright \qpp \right) q^{-1/2}
\omega_{-}, 
$$
leading to $(\dd \qpp\wedge
\dd \,\qpp)\qpp = (\del \qpp\wedge \delb \,\qpp)\, \qpp=
(E\triangleright\qpp) (F\triangleright\qpp)\,
q^{-2}\omega_{+}\wedge\omega_{-}$ again, and similarly for the case $n\leq0$.
\end{proof}
\end{lemm}

\begin{prop}
The curvature of the connection on $\ce_n=(\Asq)^{\mn+1}\qpp$ for  the canonical left covariant 2D
calculus on $\Asq$,
is given by
\beq\label{curv}
F_\nabla =  q^{-n-1} [n]\, \qpp\, \omega_{+}\wedge\omega_{-}.
\eeq
Moreover, with the left $\Asq$-modules isomorphism $\cl_n \simeq
\ce_n$ of Prop.\ref{isoeqsec} extended in a
natural way to forms $\Omega(\sq) \otimes_{\Asq} \cl_n \simeq \Omega(\sq)
\otimes_{\Asq} \ce_n$, the curvature on $\cl_n$
corresponds to the curvature on $\ce_n$, that is,
$$
F_\nabla \sigma_{\phi}  = \sigma_{\nabla^2\phi} ,
$$
for $\phi\in\cl_n$ and corresponding section $\sigma_{\phi}\in\ce_n$.
\begin{proof}
The second statement is a direct consequence of the first one in \eqref{curv} and the latter is evident once one substitute the result of Lemma~\ref{lemproj} in the expression \eqref{curvsec}. 
\end{proof}
\end{prop}

\section{The winding numbers}\label{se:wn}

The line bundles on the sphere $\sq$ described in Sect.~\ref{se:avb} are classified by their winding number $n\in\IZ$. In this section we first recall how to compute this number 
by means of a Fredholm module for the sphere. On the other hand, in order to integrate the gauge curvature on the quantum sphere $\sq$ one needs a `twisted integral'; the result is not an integer any longer but rather its q-analogue. 

\bigskip

The projections $\qpp$ given in Sect.~\ref{se:avb}, which describe the line bundles, are representatives of classes in the $K$-theory of $\sq$, i.e. $[\qpp]\in K_0(\sq)$. A way to compute the corresponding winding number is by pairing them with a nontrivial element in the dual $K$-homology, i.e. with (the class of) a nontrivial Fredholm module $[\mu]\in  K^0(\sq)$.
In fact,  it is more convenient to first compute the corresponding Chern characters in the cyclic homology $\chern_*(p) \in \mathrm{HC}_*(\sq)$ and cyclic cohomology $\chern^*(\mu)\in \mathrm{HC}^*(\sq)$ respectively, and then use the pairing between cyclic homology and cohomology.

The Chern character of the projections $\qpp$ has a non trivial component in degree zero $\chern_0(\qpp)\in \mathrm{HC}_0(\sq)$ simply given by a (partial) matrix trace:
$$
\chern_0(\qpp) := \tr (\qpp) = 
\begin{cases}
\sum_{\mu=0}^{n} \beta_{n,\mu}(c^*c)^\mu 
\prod_{j=0}^{n-\mu-1} (1-q^{-2j} c^*c), \qquad n\geq 0 \\
~ \\
\sum_{\mu=0}^{\mn} \alpha_{n,\mu}(c^*c)^{\mn-\mu} 
\prod_{j=0}^{\mu-1} (1-q^{2j} c^*c), \qquad n\leq 0 
\end{cases}  ,
$$
and $\chern_0(\qpp)\in\Asq$. Dually, one needs a cyclic 0-cocycle, that is a trace on $\Asq$. This was obtained in \cite{MNW} and it is a trace on $\Asq / \IC$, that is it vanishes on $\IC \subset \Asq$. On the other hand, its value on powers of the element $(c^*c)$ is:
$$
\mu\left((c^*c)^k\right) = (1-q^{2k})^{-1}, \qquad k>0 .
$$
The pairing was computed in \cite{H00} and results in  
$$
\hs{[\mu]}{[\qpp]}:= \mu\left(\chern_0(\qpp)\right)= -n.
$$
The integer above is a topological quantity that depends only on the bundle, both on the quantum sphere and on its classical limit, an ordinary 2-sphere. 
In this limit it could also be computed by integrating the curvature 2-form of a connection (indeed any connection) on the sphere. As mentioned, in order to integrate the gauge curvature on the quantum sphere $\sq$ one needs a `twisted integral'; furthermore the result is not an integer any longer but rather a q-integer. To proceed we need some additional ingredients.

Given a $*$-algebra $\ca$ with a state $\varphi$, 
an automorphism $\vartheta$ of $\ca$  is called a modular automorphism associated with $\varphi$ if it happens that 
$
\varphi(f g) = \varphi(\vartheta(g)f),
$ 
for $f,g\in\ca$. It is known \cite[Prop.~4.15]{KS97}, that the modular automorphism associated with the Haar state $h$ on the algebra $\ASU$ is:
\beq\label{masu2}
\vartheta(g) = K^{-2} \lt g \rt K^2 .
\eeq
The restriction of $h$ to $\Asq$ yields a faithful, invariant --
that is $h(f\rt X)=h(f) \varepsilon(X)$ for $f\in\Asq$ and $X\in\su$ --, state on $\Asq$ with modular automorphism 
\beq\label{masp}
\vartheta(g) = g \rt K^2, \qquad \mathrm{for} \quad g\in\Asq,
\eeq
the restriction of \eqref{masu2} to $\Asq$. It was proven in \cite{SW04} that, with $\omega_{+}\wedge\omega_{-}$ the central generator of $\Omega^2({\sq})$, $h$ the Haar state on $\Asq$ and $\vartheta$ its modular automorphism in \eqref{masp}, the linear functional 
\beq\label{int}
\int : \;\; \Omega^2({\sq}) \to \IC, \qquad \int f\, \omega_{+}\wedge\omega_{-} := h(f), 
\eeq
defines a non-trivial $\vartheta$-twisted cyclic $2$-cocycle $\tau$ on $\Asq$:
\beq\label{twcc}
\tau(f_0,f_1,f_2):=\int f_0\, \dd f_1 \wedge \dd f_2 .
\eeq
That is $b_\vartheta \tau=0$ and $\lambda_\vartheta \tau=\tau$ where $b_\vartheta$ is the $\vartheta$-twisted coboundary operator:                  
$$
(b_\vartheta\tau)(f_0,f_1,f_2,f_3) := \tau(f_0f_1,f_2,f_3)-\tau(f_0,f_1f_2,f_3)+\tau(f_0,f_1,f_2f_3)-\tau(\vartheta(f_3)f_0,f_1,f_2) ,
$$
and $\lambda_\vartheta$ is the $\vartheta$-twisted cyclicity operator:
$$
(\lambda_\vartheta\tau)(f_0,f_1,f_2) := \lambda_\tau(\vartheta(f_2),f_0,f_1) .
$$
The non-triviality means that there is no 1-cochain $\alpha$ on $\Asq$ such that
$b_\vartheta \alpha = \tau$ and $\lambda_\vartheta \alpha= \alpha$. Here the operators $b_\vartheta$ and $\lambda_\vartheta$ are defined by formul{ae} like the above (and directly generalizable in any degree). Thus $\tau$ is a class in $\mathrm{HC}^2_\vartheta(\sq)$, the degree 2 twisted cyclic cohomology of the sphere $\sq$.

To couple the twisted cocycle $\tau$ with the bundles over  $\sq$, one needs a twisted Chern character. In fact, for our specific case we do not need the full theory and it is enough to consider the lowest term, that of a twisted or `quantum trace' \cite{wa07}. If $M \in \Mat_{m+1}(\Asq)$, its (partial) quantum trace is the element $\qtr(M) \in \Asq$ given by
$$
%\label{qtrace}
\qtr(M):= \tr\left(M \sigma_{m/2} (K^{2})\right) :=\sum\nolimits_{jl} M_{jl} 
\left(\sigma_{m/2}(K^{2})\right)_{lj},
$$
where $\sigma_{m/2} (K^{2})$ is the matrix from \eqref{eq:uqsu2-repns}
for the spin $J=m/2$ representation of $\su$.
The q-trace is `twisted' by the automorphism $\vartheta$, that is 
$$
\qtr(M_1 M_2 ) = \qtr\left( (M_2 \rt K^{2}) M_1 \right) = 
\qtr\left( \vartheta(M_2) M_1 \right).
$$
For this one uses `right crossed product' rules:
$$
x h = {\co{h}{1}}\,(x \rt \co{h}{2}), \qquad\qquad \mathrm{for} \quad x\in\Asq, 
\quad h\in \su. 
$$
Then,
\begin{align*}
\qtr(M_1 M_2 ) &= \tr\left( M_1 M_2 \sigma_{m/2} (K^{2})\right) 
= \tr\left( M_1 \sigma_{m/2} (K^{2}) (M_2 \rt K^{2}) \right) \\
&= \tr\left((M_2 \rt K^{2}) M_1 \sigma_{m/2} (K^{2})  \right) =
\qtr\left( (M_2 \rt K^{2}) M_1 \right) .
\end{align*}

\begin{lemm}\label{qtr}
Let $\qpp$ be the projection given in \eqref{qpro} and \eqref{qpron} for $n\geq0$ or
$n\leq0$,  respectively. Then 
$$
\qtr \qpp\ = q^{n} .
$$
\begin{proof}
We prove this for the case $n\geq0$, the other case being similar. 
The matrix $\sigma_{{n}/{2}} (K^{2})$ from the expression \eqref{eq:uqsu2-repns} for the spin $J=n/2$ representation is diagonal with entries $(\sigma_{{n}/{2}} (K^{2}))_{\mu\mu}=q^{n-2\mu}$ for $\mu=0, \dots, n$. One computes
\begin{align*}
\qtr \qpp 
&= q^n 
\sum\nolimits_{\mu=0}^{n} 
q^{-2\mu} \beta_{n,\mu}\, c^{* \mu} a^{* n-\mu} a^{n-\mu} c^{\mu} 
= q^n
\sum\nolimits_{\mu=0}^{n} 
q^{-2\mu + 2\mu(n-\mu ) } \beta_{n,\mu}\,  a^{* n-\mu} c^{* \mu}  c^{\mu}
a^{n-\mu} \\
& = 
q^n \sum\nolimits_{\mu=0}^{n} 
\alpha_{n,n-\mu}\,  a^{* n-\mu} c^{* \mu}  c^{\mu}
a^{n-\mu} =
 q^n (a^{*} a + c^{*} c)^{n} = q^n, 
\end{align*}
having used the relation \eqref{remref} for the coefficients $\alpha$'s and $\beta$'s, and the identity \eqref{id2}. 
\end{proof}
\end{lemm}
We are ready to integrate the gauge curvature.
\begin{prop}
Let $F_\nabla$ be the curvature of the connection on $\ce_n=(\Asq)^{\mn+1}\qpp$ for  the canonical left covariant 2D
calculus on $\Asq$. Then, for its integral one finds:
\beq\label{intcurv}
-q \int \qtr (F_\nabla)  = - [n]\, .
\eeq
\begin{proof}
{}From eq. \eqref{curv} and Lemma~\ref{qtr} one has $\qtr(F_\nabla)=q^{-n-1} [n] \qtr(\qpp)  \omega_+\wedge\omega_- = q^{-1}[n]  \omega_+\wedge\omega_-$ and the result follows from the definition of the integral in \eqref{int} and the fact that $h(1)=1$ for the Haar state  $h$ on $\Asq$
\end{proof}
\end{prop}

\begin{rema}
{}From the definition \eqref{twcc} of the $\vartheta$-twisted cyclic $2$-cocycle $\tau$ and the expression \eqref{curv} of the curvature $F_\nabla$, the \eqref{intcurv} is also the coupling of the cocycle $\tau$ with the projection $\qpp$:
$$
(-q\, \tau) \circ  \qtr (\qpp,\qpp,\qpp) = -[n].
$$ 
In \cite{NT,wa07} this was obtained as the $\q-in$ of the Dirac operator on the sphere $\sq$. 
\end{rema}

\section{The Gauged Laplacian operator on the sphere $\sq$}
\label{se:glo}

Generalized Laplacian operators on the quantum sphere $\sq$ were
already studied in \cite{pod89}. A Hodge $\star$-operator entered
the game in \cite{maj05}.  We first recall the natural scalar
Laplacian on $\sq$ before going to its gauged version.

\subsection{A Laplacian on the quantum $\sq$ sphere}

For the Hodge $\star$-operator on 1-forms one needs a left-covariant 
map $\star:\Omega^{1}(\sq) \to \Omega^{1}(\sq)$ whose
square is the identity. In the
description of the calculus as in Prop.\ref{2dsph}, one
defines: 
\beq\label{hod1} 
\star (\del f)=\del f , \qquad \star
(\delb f)=-\delb f , \qquad \forall \, f\in \Asq , 
\eeq 
and shows its compatibility with the bimodule structure on forms.
Thus $\star$ is taken to have values $\pm 1$ on holomorphic or anti-holomorphic 1-forms respectively:  $\star\Omega^{1}_{\pm}(\sq) = \pm \Omega^{1}_{\pm}(\sq)$; 
in particular $\star \omega_\pm = \pm \omega_\pm$.
The
calculus has one central top 2-form and the above operator is naturally extended by requiring that
\beq\label{hod2} 
\star 1=\omega_{+}\wedge\omega_{-}, \qquad \star
\left(\omega_{+}\wedge\omega_{-}\right)=1 . \eeq We have all the
ingredients to define a (scalar) Laplacian operator on $\sq$:
$$
\Box f\,:=\,- \half\star \dd \star \dd f
=\,-\delb\del f=\del\delb f ,  \qquad \forall \, f\in \Asq. 
$$
Simple manipulations  yield: 
$$
\Box f=
\half\left(X_{+}X_{-}+q^{-2}X_{-}X_{+}\right) \triangleright f ,
\qquad \forall \, f\in \Asq. 
$$ 
Indeed, with
$f\in\Asq$,  one has:
$$
\dd f =\left(X_{+}\triangleright
f\right)\omega_{+}+\left(X_{-}\triangleright
f\right)\omega_{-}, \quad \mathrm{and} \quad
\star(\dd f) =\left(X_{+}\triangleright f\right)\omega_{+}-\left(X_{-}\triangleright f\right)\omega_{-};
$$
then, using the expression \eqref{d1f} for the exterior derivative
of a 1-form:
$$
\dd(\star \dd f) = - \left(X_{+}X_{-}\triangleright
f+ q^{-2}X_{-}X_{+}\triangleright f\right)\omega_{+}\wedge\omega_{-},
$$
and finally $\star(\dd\star \dd f) = -
\left(X_{+}X_{-}\triangleright f + q^{-2}X_{-}X_{+}\triangleright f
\right)$. It comes as no surprise that the (scalar) Laplacian is the
quadratic Casimir of $\su$ given in \eqref{cas}; more precisely,
\beq\label{lapcas} 
q\, \Box = C_q + \tfrac{1}{4} - [\half]^2 
. \eeq 
Indeed, since $K \lt$ acts as the
identity on elements on $\Asq$, on this algebra one has that
$\left(X_{+}X_{-}+q^{-2}X_{-}X_{+}\right)\lt = q^{-1} \left(EF +
FE\right) K^2 \lt = q^{-1} \left(EF + FE\right) \lt = 2 q^{-1}
(EF)\lt$, from the commutation relations \eqref{relsu}. The claimed
identity follows from direct comparison with the Casimir operator
\eqref{cas} restricted to $\Asq$. Relation \eqref{lapcas} makes it
also clear that the above Laplacian is related to the square of a
Dirac operator \cite{SW04,maj05}.

The spectrum and the eigenspace decomposition are computed using the
decomposition \eqref{decsp} for $\Asq=\cl_0=\bigoplus_{J\in\IN}
V_{J}^{(0)}$, for the right action of $\su$. Since left and right
action commute, this action clearly leaves
invariant the eigenspaces of the Laplacian: $\Box (f\rt  h) = (\Box f)\rt  h$, for all $h\in \su$.
We know from \eqref{hwv} that for fixed $J$ the highest weight vector
is $c^Ja^{*J}$ on which a direct computation gives:
$$
\Box(c^Ja^{*J})=q^{-1}[J][J+1]\, (c^Ja^{*J}) .
$$
The $(2J+1)$ corresponding eigenfunctions are obtained with the
action of the raising operator $\rt  F$ and are given by
$\{ (c^Ja^{*J})\rt  F^l, \, l=0,1,\dots,2J \}$. This
result is of course consistent with the equality \eqref{lapcas}.

\subsection{The gauged Laplacian operator}

\begin{defi}
Let $\nabla:\ce \to \Omega(\sq) \otimes_{\Asq} \ce$ be a covariant
derivative on the module $\ce$, with $\Omega(\sq)$ the left covariant calculus on the sphere $\sq$ described in Sect.~\ref{se:cals2}.
And let  $\star$ be the Hodge operator on $\Omega(\sq)$ as given in \eqref{hod1} and
\eqref{hod2}. The gauged Laplacian operator $\Box_\nabla : \ce \to
\ce$ is defined as:
$$
\Box_\nabla := - \half \star \nabla \star \nabla.
$$
\end{defi}

We shall presently give it explicitly on the line bundle $\cl_n$ and
on the corresponding sections $\ce_n:=(\Asq)^{\mn+1}\qpp$ showing
that they correspond to each other.
\begin{prop}
Let $\ce$ be the line bundle $\cl_n$ defined in \eqref{libu}, with
connection $\nabla$ given in \eqref{coder2d}. Then: 
\beq\label{qlpo}
\Box_\nabla \phi = \half q^{-2n}
\left(X_{+}X_{-}+q^{-2}X_{-}X_{+}\right)\triangleright\phi,
\qquad\mathrm{for} \quad \phi\in\cl_n. 
\eeq 
Moreover, with the left
$\Asq$-modules isomorphism $\cl_n \simeq \ce_n:=(\Asq)^{\mn+1}\qpp$
of Prop.\ref{isoeqsec}, the Laplacian on $\cl_n$ corresponds to the
one on $\ce_n$ for the Grassman connection, that is:
$$
\Box_\nabla (\sigma_{\phi} ) = \sigma_{\Box_\nabla\phi},
$$
for $\phi\in\cl_n$ and corresponding section $\sigma_{\phi}\in\ce_n$.
\begin{proof}
Using \eqref{coder2d} and \eqref{hod1} we have, $\star \nabla \phi =
q^{-n-2} \omega_{+} \left(X_{+}\triangleright\phi\right) - q^{-n+2}
\omega_{-} \left(X_{-}\triangleright\phi\right)$ and (being
$\dd\omega_-=\dd\omega_+=0$ on $\sq$),
\begin{align*}
\nabla(\star \nabla \phi) &= - q^{-n-2} \omega_{+} \wedge \left(X_{-} X_{+}\triangleright\phi\right) \omega_-
+q^{-n+2} \omega_{-} \wedge \left(X_{+} X_{-}\triangleright\phi\right))\omega_{+} \\
& = - q^{-2n-2} \omega_{+}\wedge\omega_- \left(X_{-} X_{+}\triangleright\phi\right)
+q^{-2n+2} \omega_{-} \wedge \omega_{+} \left(
X_{+} X_{-}\triangleright\phi\right) \\
& = - q^{-2n} \omega_{+}\wedge \omega_-  \left[ q^{-2} \left(X_{-} X_{+}\triangleright\phi\right) +
 \left(X_{+} X_{-}\triangleright\phi\right) \right] ,
\end{align*}
and using \eqref{hod1} the first statement follows.

For the second point, if $\sigma\in\ce_n$ and $\nabla
\sigma=\dd(\sigma) \qpp$ is the Grassmann connection, a computation
similar to the previous one leads to:
\begin{multline}\label{12}
\Box_\nabla \sigma =
\half
\left[\left( X_{+}X_{-} +q^{-2}X_{-}X_{+} \right) \lt \sigma \right] \qpp \\ +
\half \left[ q^{-2} (X_{+}\lt \sigma) X_{-}\lt \qpp - (X_{-}\lt \sigma) X_{+}\lt \qpp \right]\qpp . 
\end{multline}
To continue, consider the case $n\geq0$, the case $n\leq0$ being similar.
{}From Prop.\ref{isoeqsec}, to $\phi\in\cl_n$ there corresponds the
section $\sigma_\phi= \phi \bra{\Psi^{(n)}} \in  \ce_n$, where
$\bra{\Psi^{(n)}}$ and its dual $\ket{\Psi^{(n)}}$ are the vector
valued functions in \eqref{qpro} with projection 
$\qpp=\qpp=\ket{\Psi^{(n)}}\bra{\Psi^{(n)}}$. Using the vanishing 
$X_{-} \lt \bra{\Psi^{(n)}}=X_{+} \lt \ket{\Psi^{(n)}}=0$ and that 
$K \lt \bra{\Psi^{(n)}}=q^{-n/2}\bra{\Psi^{(n)}}$ and $K \lt \ket{\Psi^{(n)}}=
q^{n/2}\ket{\Psi^{(n)}}$, we compute:
\begin{multline*}
\half \left(X_{+}X_{-} + q^{-2}X_{-}X_{+}\right) \lt \sigma_\phi  = \half
q^{-2n} \left[ \left(X_{+}X_{-} + q^{-2}X_{-}X_{+}\right) \lt \phi \right] \bra{\Psi^{(n)}} \\
+ q^{-n} (X_{-} \lt \phi) \, X_{+} \lt \bra{\Psi^{(n)}} + \half
q^{-2}  \phi \, (X_{-} X_{+} \lt \bra{\Psi^{(n)}} ),
\end{multline*}
and using  $\left( X_{+} \lt \bra{\Psi^{(n)}} \right)
\ket{\Psi^{(n)}}=0$, we arrive at:
\begin{multline}
\half \left[ \left(X_{+}X_{-} + q^{-2}X_{-}X_{+}\right) \lt \sigma_\phi\right] \qpp  \label{13} \\
\qquad\qquad = \half
q^{-2n} \left[ \left(X_{+}X_{-} + q^{-2}X_{-}X_{+}\right) \lt \phi \right] \bra{\Psi^{(n)}}
+ \half
q^{-2}  \phi \, (X_{-} X_{+} \lt \bra{\Psi^{(n)}} ) \qpp . 
\end{multline}
On the other hand, we know that $(X_{+}\lt \qpp )
\qpp=0$ and that $(X_{-}\lt \qpp ) \qpp = X_{-}\lt \qpp = q^{-n}
\left( X_{-} \lt \ket {\Psi^{(n)}} \right) \bra{\Psi^{(n)}} $. Using
 $\hs{\Psi^{(n)}} {X_{-} \lt {\Psi^{(n)}}}=0$, the
second part in \eqref{12} becomes:
$$
-\half q^{-n-2} \phi \, \left(X_{+}\lt
\bra{\Psi^{(n)}} \right) \left( X_{-}\lt \ket{\Psi^{(n)}} \right) \bra{\Psi^{(n)}},
$$
and this cancel the second term in \eqref{13} due to:
\begin{align*}
& \left(X_{-} X_{+} \lt \bra{\Psi^{(n)}} \right) \ket{\Psi^{(n)}} +
q^{-n} \left(X_{+}\lt \bra{\Psi^{(n)}} \right) \left(X_{-}\lt \ket{\Psi^{(n)}} \right) \\
& \qquad = q^{-n} \left\{ X_{-}\lt\left(X_{+} \lt \bra{\Psi^{(n)}} \right) \left(K^2\lt \ket{\Psi^{(n)}}\right) +
\left(X_{+}\lt \bra{\Psi^{(n)}} \right) \left(X_{-}\lt \ket{\Psi^{(n)}} \right) \right\} \\
& \qquad = q^{-n} X_{-} \lt \left[ \left( X_{+} \lt \bra{\Psi^{(n)}} \right) \ket{\Psi^{(n)}}\right] = 0.
\end{align*}
Collecting all the above we arrive at:
$$
\Box_\nabla (\sigma_\phi) = \half
q^{-2n} \left[ \left(X_{+}X_{-} + q^{-2}X_{-}X_{+}\right) \lt \phi \right] \bra{\Psi^{(n)}} = \sigma_{\Box_\nabla\phi}
$$
and this ends the proof.
\end{proof}
\end{prop}

We are ready to diagonalize the gauged Laplacian $\Box_\nabla$. For this we use the decomposition \eqref{decoln} of the modules $\cl_n$ for the right action of $\su$, as this action leaves invariant the eigenspaces of the
gauged Laplacian: $\Box_\nabla (\phi \rt  h) = (\Box_{\nabla} \phi)\rt  h$,
for all $h\in \su$. 

{}From Sect.~\ref{qdct} we know that the integer $n$ (the monopole charge) labels topological sectors and given $n$ the admissible values of $J$ in the decomposition $\cl_n=\bigoplus V_{J}^{\left(n\right)}$ are 
$J=\frac{\mn}{2}+s$, 
with $s\in\IN $. Moreover, the highest weight elements in each $V_{J}^{\left(n\right)}$ is given in \eqref{hwv} as  
$\phi_{n,J}=c^{J-n/2}a^{*J+n/2}$, and the remaining $2J$
basis vectors in $V_{J}^{\left(n\right)}$ are obtained via the right action
$\triangleleft E$. The vectors
$\phi_{n,J,l}=(c^{J-n/2} a^{*J+n/2})\rt  E^{l}$, with $l=0,\ldots,
2J$, are eigenfunctions of $\Box_\nabla$ with the same eigenvalue. 
\begin{prop}
On the vectors $\phi_{n,J,l}$ the gauged Laplacian is diagonal,
$$
\Box_\nabla\phi_{n,J,l} = \lambda_{n,J} \phi_{n,J,l},
$$ 
for $l=0,\ldots, 2J$, with the $(2J+1)$-degenerate energies:
\begin{align}\label{energies}
\lambda_{n,J} 
& = q^{-n-1} \left( [J-\tfrac{n}{2}]\, [J+\tfrac{n}{2}+1] + \half [n] \right) \nn \\
& = q^{-n-1} \left\{ [J+\half]^2 - \half \left( [\tfrac{n+1}{2}]^2 + [\tfrac{n-1}{2}]^2 \right) \right\} .
\end{align}
\begin{proof}
A short computation yields:
\begin{align*}
\Box_\nabla (c^{J-n/2} a^{*J+n/2})  &= 
\half q^{-2n} \left(X_{+}X_{-}+q^{-2}X_{-}X_{+}\right)\triangleright (c^{J-n/2} a^{*J+n/2}) \\
&=\lambda_{n,J} (c^{J-n/2}a^{*J+n/2}), 
\end{align*}
with $\lambda_{n,J}$ given in the first line of \eqref{energies}. The second equality there is obtained with a direct algebraic manipulation. 
\end{proof}
\end{prop}

A remarkable fact is that, contrary to what happens in the classical limit, the energies are not symmetric under the exchange $n \leftrightarrow -n$, an additional example that `quantization removes degeneracy'. 
Writing $J=\frac{\mn}{2}+s$, with $s\in\IN$, the energies become:
$$
\lambda_{n,s}=q^{-n-1}\left([s][n+s+1] + \half [n]\right) , 
\qquad \mathrm{for} \quad n\geq0 ,
$$
with $(n+2s+1)$ eigenfunctions $\phi_{n,s,l}=(c^{s} a^{*n+s})\rt E^l$, and 
$$
\lambda_{n,s}=q^{-n-1}\left([s-n][s+1] + \half[n] \right) ,
\qquad \mathrm{for} \quad n\leq0 ,
$$
with $(\mn+2s+1)$ eigenfunctions $\phi_{n,s,l}=(c^{s+\mn}a^{*s})\rt E^l$.

Having in mind a physics
parallel with the quantum Hall effect, the integer $s$ labels Landau
levels and the $\phi_{n,s,l}$ are the (`one excitation') Laughlin wave functions with
energies $\lambda_{n,s}$. The lowest Landau, $s=0$,  
is $\mn$-degenerate with energy $\lambda_{n,0} = \half q^{-n-1}
[\mn]$.

It is worth spending a few words on the classical limit. At the value $q=1$, the energies of the gauged Laplacian become
\begin{align*}
%\label{energies0}
\lambda_{n,s}(q\to1) & =  (J-\tfrac{n}{2})(J+\tfrac{n}{2}+1) + \half n 
= J(J+1) - \tfrac{1}{4} n^2 \nn \\
& = (\half \mn + s)(\half \mn + s +1) - \tfrac{1}{4} n^2 = \mn ( s+ \half) + s(s+1), 
\end{align*}
and coincide with the energies of the classical gauged Laplacian (see e.g. \cite{H83}). Clearly, they are symmetric under the exchange  $n \leftrightarrow -n$ which corresponds to inverting  the direction of the magnetic field.

\bigskip
The relation of the gauged Laplacian with the quadratic Casimir of $\su$ is more involved than its classical counterpart.
\begin{prop}\label{lapcasn}
When acting on the module $\cl_n$ of equivariant maps the gauged Laplacian and the quadratic Casimir of $\su$ are related as:
\beq\label{lc1}
q^{n+1}\Box_{\nabla} = C_{q} + \tfrac{1}{4} -
\tfrac{1}{2}\left(\left[\tfrac{n+1}{2}\right]^{2}+\left[\tfrac{n-1}{2}\right]^{2}\right)\,.
\eeq
Moreover, 
the operators $\Box_{\nabla}$ and $C_{q}$ are related as:
\beq\label{lc2}
q K^2 \Box_{\nabla} = C_{q} + \tfrac{1}{4} - \half 
\left(\frac{qK^{2}-2+q^{-1}K^{-2}}{(q-q^{-1})^{2}} +
\frac{q^{-1}K^{2}-2+qK^{-2}}{(q-q^{-1})^{2}}\right) 
\, .
\eeq
\begin{proof}
We first note that the operator $K$ commutes with all operators involved in the proposition and we can have it either on the left or on the right. Then, in \eqref{qlpo} giving the action of the gauged Laplacian on an equivariant map, 
the factor $q^{-2n}$ can be traded for the action of the operator $K^{-4}$. In turn, 
$$
\Box_{\nabla} = 
\half K^{-4}
(X_{+}X_{-}+q^{-2}X_{-}X_{+}) 
= \half q^{-1}K^{-4}(EF+FE)K^{2} 
= \half q^{-1} K^{-2} (EF+FE) .
$$
Denote $\lambda=(q-q^{-1})^{-1}$ for simplicity. For the 
Casimir operator in \eqref{cas} we get: 
\begin{align*}
C_{q} 
&= \half(FE+EF) + \half [F,E] + \lambda^2 (qK^{2}-2+q^{-1}K^{-2})-\tfrac{1}{4} 
\\
&= \half(FE+EF) + \lambda^2 (qK^{2}-2+q^{-1}K^{-2}) 
- \half \lambda (K^{2}-K^{-2}) -\tfrac{1}{4} \\
& = \half(FE+EF) + \half \lambda^2 (qK^{2}-2+q^{-1}K^{-2} + 
q^{-1}K^{2}-2+qK^{-2}) -\tfrac{1}{4}
\,.
\end{align*}
A comparison between these two expressions establishes \eqref{lc2} which in turn, when computed on $\cl_n$ gives \eqref{lc1}.
\end{proof}
\end{prop}

The above proposition has the expected classical limit. Again, 
by setting $K=q^H$, expanding in the parameter $\hbar=:\log q$ and truncating at the 0-th order in $\hbar$, from a direct computation the relation \eqref{lc2} becomes
$$
\Box_\nabla = C_{\mathrm{SU(2)}} - H^2. 
$$
This is an example of the general result recalled in the introduction, being $H$ the generator of the structure group $\U(1)$ with Casimir element $C_{\U(1)}=H^2$.

%\section{Final remarks}

\subsection*{Acknowledgements}
This work was partially supported by the `Italian project Cofin06 - Noncommutative geometry, quantum groups and applications'.
The research of AZ started at SISSA (Trieste, Italy) and went on at the IAM at
Bonn University (Germany), thanks to a fellowship by the Alexander von
Humboldt Stiftung; he thanks the Mathematical Physics Sector
of SISSA and his host in Germany, Prof. Sergio Albeverio, for their warm
hospitality. We thank Francesco D'Andrea for reading the comptu-script.
%\newpage
\appendix

\section{An appendix of preliminaries}\label{se:pre}
We need to recall few facts about quantum principal bundles with (nonuniversal)
differential calculi and connections on them. We briefly review their main properties here starting 
with covariant differential calculi.

\subsection{Differential calculi}

Given a $\IC$-algebra with unit $\ca$, any first order
differential calculus $(\oca{1}, \dd)$ on $\ca$ can be obtained from
the universal calculus $(\oca{1}_{un}, \delta)$. The space of universal 1-forms
is the submodule of $\ca \otimes \ca$ given by
$\oca{1}_{un} := \ker(m:\ca \otimes \ca \to \ca)$, with $m(a\otimes
b)=ab$ the multiplication map. The universal differential 
$\delta: \ca \to \oca{1}_{un}$ is $\delta a := 1 \otimes a - a \otimes 1$.
If $\cn$ is any sub-bimodule of $\oca{1}_{un}$ with
projection $\pi: \oca{1}_{un} \to \oca{1} := \oca{1}_{un} /\cn$, then $(\oca{1}, \dd)$, with $\dd:=\pi\circ\delta$, is a first order
differential calculus over $\ca$ and any such a calculus can be
obtained in this way.

If the algebra $\ca$ is covariant for the coaction of a quantum group $\ch=(\ch,\Delta, \eps, S)$, one has a notion of covariant calculi on $\ca$ 
as well. 
Then, let  $\ca$ be a (right, say) $\ch$-comodule algebra, with a right  coaction $\Delta_R : \ca \to \ca \otimes \ch$ which is also an algebra map. In order to state the covariance of the calculus 
$(\oca{1}, \dd)$ one needs to extend the coaction of $\ch$. A map 
$  \Delta_{R}^{(1)} : \oca{1} \to \oca{1}\otimes \ch$ is defined by the requirement
$$
\Delta_{R}^{(1)}(\dd f) = (\dd \otimes \id) \Delta_{R}(f)
$$
and bimodule structure governed by
$$
\Delta_{R}^{(1)}(f \dd f') = \Delta_{R}(f) \Delta_{R}^{(1)}(\dd f') , \qquad 
  \Delta_{R}^{(1)}((\dd f) f') = \Delta_{R}^{(1)}(\dd f) \Delta_{R}(f').
$$
The calculus is said to be right covariant it  happens that 
$$
(\id \otimes \Delta) \Delta_{R}^{(1)} = (\Delta_{R}^{(1)} \otimes \id) \Delta_{R}^{(1)} \qquad \mathrm{and} \qquad (\id \otimes \eps) \Delta_{R}^{(1)} = 1 .
$$
A calculus is right covariant if and only if for the corresponding bimodule $\cn$ it is verified that $\Delta_{R}^{(1)}(\cn) \subset \cn \otimes \ch$, where $\Delta_{R}^{(1)}$ is defined on  $\cn$ by formul{\ae} as above with the universal derivation $\delta$ replacing the derivation $\dd$. Right covariance of the calculus implies that $\oca{1}$ has a module basis $\{\eta_{a}\}$ of right invariant 1-forms, that is 1-forms for which
$$
\Delta_{R}^{(1)}(\eta_{a})= \eta_{a} \otimes 1. 
$$
Similar consideration and formul{\ae} holds for left covariance under a left coaction $\Delta_L$. 
In particular, left covariance of a calculus similarly implies that
$\oca{1}$ has a module basis $\{\omega_{a}\}$ of left invariant 1-forms, that is 1-forms for which, 
$\Delta_{L}^{(1)}(\omega_{a})=1\otimes\omega_{a}$.

Differential calculi on a quantum group $\ch=(\ch,\Delta, \eps, S)$ were already studied in \cite{wor89}. Now $\Delta: \ch \to \ch \otimes \ch$ is viewed as both a right and a left coaction of $\ch$ on itself. Right and left covariant calculi on $\ch$ will be defined as before with in particular a basis of invariant forms for the corresponding covariant calculus. In addition one has the notion of a bicovariant (that is both right and left covariant) calculus. 
%
%Let be a quantum group (a Hopf algebra). If
%$(\oca{1}, \dd)$ is a first order differential calculus on $\ch$,
%the left coaction of $\ch$ is extended to $\oca{1}$ as
%\beq\label{lcov} \Phi_{L}(f\dd f^{\prime}) =\Delta(f)(\id \otimes
%\dd)\Delta(f^{\prime}) . \eeq The given calculus is left covariant
%when this left coaction is compatible with the coproduct and with
%the counit in $\ch$, and the coproduct extends to a left coaction on
%$\oca{1}$ such that $\dd$ is an intertwiner and $\Phi_{L}$ is a
%bimodule map. Left covariance of a calculus also implies that
%$\oca{1}$ has a module basis of left invariant 1-forms, that is 1
%forms for which, \beq \Phi_{L}(\omega_{a})=1\otimes\omega_{a} . \eeq

Not surprisingly, on a quantum group there is more structure. 
Given the bijection 
$$
r : \ch\otimes\ch\mapsto\ch\otimes\ch,
\qquad r(h\otimes h^{\prime}):=(h\otimes1)\Delta(h^{\prime}) , 
$$ 
one proves that $ r(\och{1}_{un})=\ch\otimes \ker\eps $. 
Then, if
$\mathcal{Q} \subset \ker\eps$ is a right ideal of $\ker\eps$, the
inverse image, $\mathcal{N}_{\mathcal{Q}}=r^{-1}(\ch\otimes\mathcal{Q})$,
is a sub-bimodule contained in $\och{1}_{un}$.
The differential calculus defined by such a bimodule,
$\och{1}:=\och{1}_{un}/\mathcal{N}_{\mathcal{Q}}$,
is left-covariant, and any left-covariant differential calculus can
be obtained in this way. Bicovariant calculi are in one to one correspondence with such right ideals $\mathcal{Q}$ which are in addition stable under the right adjoint coaction $\Ad$ of $\ch$ onto itself, that is $\Ad(\cq) \subset \cq \otimes \ch$.  
Explicitly, one has $\Ad=\left(\id \otimes m \right) \left(\sigma\otimes
\id \right)\left(S\otimes\Delta\right)\Delta$, with $\sigma$ the flip
operator, or $\Ad(h) = \co{h}{2} \otimes \left(S(\co{h}{1}) \co{h}{3} \right)$
using the Sweedler notation $\cop h =: \co{h}{1} \otimes \co{h}{2}$ with summation understood, and higher numbers for iterated coproducts.
%
%Moreover, it is possible to prove that the dimension of the
%differential calculus is given by the dimension of the quotient
%$(Kern\,\varepsilon_{\ch})/\mathcal{Q}_{\ch}$.
%

\bigskip
The ideal $\mathcal{Q}$ also determines the quantum tangent space of
the calculus. This is a collection $\{X_{a}\}$ of elements in
$\cu_q({\ch})$ -- the Hopf algebra dual to $\ch$ -- which allows one to
write the exterior differential as 
$$
\dd h := \sum\nolimits_a
~(X_{a} \triangleright h) ~\omega_{a} , 
$$
for $h\in\ch$ and elements $X_{a}$ acting on the left on $h$. To be more
specific, with the dual pairing $\hs{~}{~}:\cu_q({\ch}) \times \ch
\to \IC$, the quantum tangent space determined by the ideal
$\mathcal{Q}$ is
$$
\mathcal{X}_{\mathcal{Q}}:=
\{X\in \ker\varepsilon_{\cu_{q}(\ch)} ~:~ \hs{X}{Q}=0, \,\, \forall \, Q\in\mathcal{Q}\},
$$
where $\varepsilon_{\cu_{q}(\ch)}$ is the counit of
$\cu_{q}(\ch)$. The left action is given by $X \lt h := (\id
\otimes X)\,\cop h$ or equivalently by $X \lt h := \co{h}{1}
\,\hs{X}{\co{h}{2}}$. The twisted
derivation nature of elements in $\mathcal{X}_{\mathcal{Q}}$ is
expressed by their coproduct, 
$
\Delta(X_{a})=1\otimes X_{a}+
\sum\nolimits_b X_{b}\otimes f_{ba}
$,
with the elements  $f_{ab} \in \cu_{q}(\ch)$ having
specific properties \cite{wor89}. These elements also control the commutation
relation between the basis 1-forms and elements of $\ch$:
$$
\omega_{a} h = \sum\nolimits_b (f_{ab}\triangleright h)\omega_{b}, \qquad
h \omega_{a} = \sum\nolimits_b \omega_{b} \left( (f_{ab}\circ S^{-1} )\triangleright h \right)  
\qquad \mathrm{for} \quad h \in \ch . 
$$

\subsection{Quantum principal bundles and connections}\label{se:qpb}

The quantum principal bundles with nonuniversal calculi we are interested in were introduced in \cite{BM93} (with refinements in \cite{BM97}). 
As a total space we consider an algebra $\mathcal{P}$ (with multiplication $m : \mathcal{P}\otimes \mathcal{P} \to \mathcal{P}$) and as structure group  
a Hopf algebra $\ch$. Thus $\mathcal{P}$ is a right $\ch$-comodule algebra
with coaction $\Delta_{R}\,:\,\mathcal{P}\mapsto\mathcal{P}\otimes\ch$.
The subalgebra of the right coinvariant elements, 
$
\mathcal{B}=\mathcal{P}^{\ch}\,:=\,\{p\in\mathcal{P}\,:\,\Delta_{R} p
= p\otimes 1\}
$,
is the base space of the bundle.
At the `topological level' the principality of the bundle is the requirement of 
exactness of the sequence:
\beq
0\, \to \,\mathcal{P}\left(\Omega^{1}(\mathcal{B})_{un}\right)\mathcal{P}\,
 \to \,\Omega^{1}(\mathcal{P})_{un}
\,\stackrel{\chi}  \to \,\mathcal{P}\otimes \ker\varepsilon_{\ch}\,  \to \,0
\label{topes}
\eeq
with $\Omega^{1}(\mathcal{P})_{un}$ and $\Omega^{1}(\mathcal{B})_{un}$ the universal calculi and the map $\chi$ is defined by, 
\beq
\chi:\,\mathcal{P}\otimes\mathcal{P}\, \to\,\mathcal{P}\otimes\ch,\qquad 
\chi:=\left(m \otimes \id\right)\left(\id\otimes\Delta_{R}\right)
\label{chimap},
\eeq
or $\chi(p' \otimes p) = p' \Delta_R(p)$. The exactness of this sequence is equivalent to the requirement that the analogous `canonical map'  $\mathcal{P}\otimes_{\mathcal{B}}\mathcal{P}\to\,
\mathcal{P}\otimes\ch$ (defined as the formula above) is an isomorphism. This is the definition that the inclusion $\mathcal{B} \hookrightarrow  \mathcal{P}$ be a Hopf-Galois extension \cite{sc90}. 
For quantum structure groups which are cosemisimple and have bijective antipodes, Th.~I of \cite{sc90} grants further nice properties. In particular, the surjectivity of the canonical map implies its bijectivity and faithfully flatness
of the extension.
%Moreover, an additional useful result \cite{scP} is that the map $\chi$ is
%surjective whenever, for any generator $h$ of $H$, the element $1
%\ot h$ is in its image.

The surjectivity of the map $\chi$ is the translation, for the deformed case, of the classical condition  that the action of the structure group on the total space of the bundle is free.

\bigskip
With differential calculi on both the total algebra $\cp$ and the structure Hopf algebra $\ch$ one needs compatibility conditions that eventually lead to an exact sequence like in \eqref{topes} with the  calculi at hand replacing the universal ones. Then, let $\left(\Omega^{1}(\mathcal{P}),\dd\right)$ be 
a $\ch$-covariant differential calculus  on $\mathcal{P}$ given via the subbimodule $\mathcal{N}_{\mathcal{P}} \in \left(\Omega^{1}(\mathcal{P})_{un}\right)$, and $\left(\Omega^{1}(\ch),\dd\right)$ a bicovariant one on $\ch$ given via the $\Ad$-invariant right ideal $\mathcal{Q}_{\ch} \in \ker\varepsilon_{\ch}$. 
The compatibility of the calculi are the requirements that 
$\chi(\mathcal{N}_{\mathcal{P}})\subseteq \mathcal{P}\otimes\mathcal{Q}_{\ch}$ and that the map 
$\sim_{ {\mathcal{N}_{{\mathcal{P}} }}} : \Omega^{1}(\mathcal{P}) \to  
\mathcal{P}\otimes(\ker\varepsilon_{\ch}/\mathcal{Q}_{\ch})$, defined by the diagram
\beq
\begin{array}{lcl}
\Omega^{1}(\mathcal{P})_{un} 
& \stackrel{{\pi_{\mathcal{N}}}} {\longrightarrow} & \Omega^{1}(\mathcal{P}) \\
\downarrow \chi  &  & \downarrow\sim_{{ \mathcal{N}_{\mathcal{P}}}} \\
\mathcal{P}\otimes \ker\varepsilon_{\ch} &
\stackrel{{\id\otimes\pi_{\cq}}}{\longrightarrow}&
\mathcal{P}\otimes(\ker\varepsilon_{\ch}/\mathcal{Q}_{\ch})
\end{array}
\label{qdia}
\eeq
(with $\pi_{\mathcal{P}}$ and $\pi_{\cq}$ the natural projections) is
surjective and has kernel 
\beq 
\ker\sim_{{ 
\mathcal{N}_{\mathcal{P}}}} = \mathcal{P}\Omega^{1}(\mathcal{B})\mathcal{P}
=:\Omega^{1}_{\mathrm{hor}}(\cp). 
\label{codd}
\eeq
Here $\Omega^{1}(\mathcal{B})=\cb \dd \cb$ is the space of nonuniversal 
1-forms on $\mathcal{B}$ associated to the bimodule $\mathcal{N}_{\mathcal{B}}:= \mathcal{N}_{\mathcal{P}} \cap \Omega^{1}(\mathcal{B})_{un}$. 
These conditions ensure the exactness of the sequence:
\beq
0\,\to\,\mathcal{P}\Omega^{1}(\mathcal{B})\mathcal{P}\,\to\,
\Omega_{1}(\mathcal{P})
\,\stackrel{\sim_{\mathcal{N}_{\mathcal{P}}}} \longrightarrow\, \mathcal{P}\otimes \left(\ker\varepsilon_{\ch}/\mathcal{Q}_{\ch}\right)\,\to\,0 .
\label{des}
\eeq
The condition 
$\chi(\mathcal{N}_{\mathcal{P}})\subseteq \mathcal{P}\otimes\mathcal{Q}_{\ch}$ is needed to have a map 
$\sim_{ {\mathcal{N}_{{\mathcal{P}} }}}$ well defined. In fact, with all conditions for a quantum principal bundle $(\cp, \cb, \ch; \cn_{\cp}, \cq_{\ch})$ satisfied, this inclusion implies equality $\chi(\mathcal{N}_{\mathcal{P}})=\mathcal{P}\otimes\mathcal{Q}_{\ch}$.
Moreover, if $(\cp, \cb, \ch)$ is a quantum principal bundle with the universal calculi, the equality $\chi(\mathcal{N}_{\mathcal{P}})=\mathcal{P}\otimes\mathcal{Q}_{\ch}$ ensures that $(\cp, \cb, \ch; \cn_{\cp}, \cq_{\ch})$ is a quantum principal bundle with the corresponding nonuniversal calculi.

%In the classical setting, the definition of a principal bundle leads to the notion of vertical vector field, namely the realisation of the Lie algebra of the gauge group in terms of vector fields on the total space manifold. In the quantum setting, this concept is in a sense dualised. 

Elements in the quantum tangent space  $\mathcal{X}_{\cq_{\ch}}(\ch)$ giving the calculus on the structure quantum group $\ch$ act on $\ker\varepsilon_{\ch} / \mathcal{Q}_{\ch}$
via the pairing $\hs{\cdot}{\cdot}$ between $U_{q}(\ch)$ and $\ch$. Then, with each $\xi\in\mathcal{X}_{\cq_{\ch}}(\ch)$ one defines a map 
\beq\label{hf}
\tilde{\xi} : \Omega^{1}(\mathcal{P}) \to\mathcal{P}, \qquad \tilde{\xi}:=\left(\id\otimes\xi\right) \circ (\sim_{\mathcal{N}_{\mathcal{P}}} )
\eeq
and declare a 1-form $\omega \in \Omega^{1}(\mathcal{P})$ to be  horizontal if{}f
$\tilde{\xi}\left(\omega\right)=0$,  for all elements $\xi\in
\mathcal{X}_{\cq_{\ch}}(\ch)$.
The collection of horizontal 1-forms is easily seen to coincide with $\Omega^{1}_{\mathrm{hor}}(\cp)$ in \eqref{codd}.

\bigskip
Finally, we come to the notion of a connection. Recall that 
the covariance condition, $\Delta_{R}\mathcal{N}_{\mathcal{P}} \subset\mathcal{N}_{\mathcal{P}} \otimes\ch $ allows one to extend the coaction $\Delta_{R}$ of $\ch$ on $\cp$ to a coaction of $\ch$ on 1-forms, 
$\Delta_{R}^{(1)} :\Omega^{1}(\mathcal{P})\mapsto\Omega^{1}(\mathcal{P})\otimes\,\ch $,
by requiring that 
$\Delta_{R}^{(1)}\circ\,\dd=(\dd\otimes\,\id)\circ \Delta_{R}$.
A connection on the 
quantum principal bundle $(\cp, \cb, \ch; \cn_{\cp}, \cq_{\ch})$ can be given as a map $\omega : \ker\varepsilon_{\ch} \to \Omega^{1}(\mathcal{P})$ such that $\xi \circ \omega = 1 \otimes \id$ and $\Delta_{R}^{(1)} \circ \omega = (\omega \otimes \id) \circ \Ad$. Here the map $\Ad$ is the quotient of the right adjoint action on $\ch$ to the space $\ker\varepsilon_{\ch}$ determined by $\Ad \circ \pi_{\cq} = (\pi_{\cq} \otimes \id) \circ \Ad$. 
It is shown in \cite{BM93,BM97} that connections are in 1-1 correspondence with $\ch$-covariant complements to the horizontal forms $\Omega^{1}_{\mathrm{hor}}(\cp)\subset\Omega^{1}(\cp)$.

\section{Computations for the quantum Hopf bundle}\label{se:cqpb}
We collect here some of the computations toward showing the quantum principal bundle structure of the $\U(1)$-bundle over the standard Podle\'s sphere $\sq$ with total space the quantum group $\SU$, given  in Sect.~\ref{se:qhb}.

\subsection{The principal bundle condition}\label{se:pbc}

We prove here that the datum $(\ASU,\Asq,\ca(\U(1)))$ is a quantum principal bundle. This is done by showing exactness of the sequence 
\begin{multline*}
0\,\to\,\ASU\left(\Omega^{1}(\Asq)_{un}\right)\ASU\, \to \\ 
\to\, \Omega_{1}(\ASU)_{un} \, \stackrel{\chi}\longrightarrow \, 
\ASU \otimes \ker\varepsilon_{\ca(\U(1))}\,\to\,0
\end{multline*}
or equivalently that the map 
$\chi : \Omega^{1}(\ASU)_{un} \to \ASU\otimes \ker\varepsilon_{\ca(\U(1))}$ defined as in (\ref{chimap}) -- and with the $\ca(\U(1))$-coaction on $\ASU$ given in (\ref{cancoa}) --, is surjective. 
Now, a generic element in $\ASU\otimes \ker\varepsilon_{\ca(\U(1))}$ is of the form $f\otimes\left(1-z^{*n}\right)$ with $n\in\IZ$ and $f\in\ASU$. To show surjectivity of $\chi$ it is enough to show that $1\otimes\left(1-z^{*n}\right)$ is in its image since left $\ASU$-linearity of $\chi$ will give the general result: if $\gamma\in\Omega^{1}(\ASU)_{un}$ is such that $\chi(\gamma)=1\otimes\left(1-z^{*n}\right)$, then $\chi(f \gamma)=f\left(1\otimes(1-z^{*n})\right)=f \otimes \left(1-z^{*n}\right)$. 

Firstly the case $n\geq0$. With $\ket{\Psi^{(n)}}$ given in \eqref{qpro}, if $\gamma\in\Omega^{1}(\ASU)_{un}$ is  
\begin{multline*}
\gamma  = \hs{\Psi^{(n)}}{\delta\Psi^{(n)}} := 
\sum\nolimits_{\mu=0}^{n}\beta_{n,\mu} a^{n-\mu} c^{\mu}
\, \delta(c^{*\mu} a^{*n-\mu}) \\
 = 1\otimes 1 - 
\sum\nolimits_{\mu=0}^{n}\beta_{n,\mu} a^{n-\mu} c^{\mu}
\, \otimes (c^{*\mu} a^{*n-\mu})
\end{multline*}
we get 
$
\chi\left(\gamma\right)\,= 1\otimes 1 - 
\sum_{\mu=0}^{n}\beta_{n,\mu} a^{n-\mu} c^{\mu}
\, (c^{*\mu} a^{*n-\mu} \otimes z^{*n}) = 1 \otimes\left(1-z^{*n}\right)
$.\\
For $n\leq0$ the proof goes analogously with the vectors $\ket{\check{\Psi}^{\left(n\right)}}$. The case $n=0$ is trivial.

\subsection{The compatibility of the calculi}\label{se:cc}

In Sect.~\ref{se:cotqpb} we have given calculi on the principal bundle
$(\ASU,\Asq,\ca(\U(1)))$. Out of the 3D left covariant calculus
$\Omega^{1}(\ASU)$, with defining ideal $\cq_{\SU}$ given in Sect.~\ref{se:lcc}, a calculus on $\ca(\U(1))$ was obtained by projection while a calculus on the subalgebra $\Asq$ was given by restriction. The `principal bundle compatibility' of these calculi is established by showing that the sequence $\eqref{des}$ is exact. For the case at hand, this sequence becomes,
\begin{multline*}
0\, \to\, \ASU\left(\Omega^{1}(\sq)\right)\ASU\, \to \\ \to\,\Omega^{1}(\ASU)\,\stackrel{\sim_{\mathcal{N}_{\SU}}}  \longrightarrow \,\ASU\otimes \ker\varepsilon_{\ca(\U(1))}/\mathcal{Q}_{\ca(\U(1))}\,\to\,0 ,
\end{multline*}
where $\mathcal{Q}_{\ca(\U(1))}$ is the ideal given in Sect.~\ref{se:csg} that defines the calculus on $\ca(\U(1))$ and the map $\sim_{\mathcal{N}_{\SU}}$ is defined as in the diagram \eqref{qdia} which now becomes,
$$
\begin{array}{lcl}
\Omega^{1}(\ASU_{un} & \stackrel{ \pi_{\cq_{\SU}} }{\longrightarrow} 
& \Omega^{1}(\ASU \\
\downarrow \chi  &  & \downarrow \sim_{\cn_{\SU}} \\
\ASU \otimes \ker\varepsilon_{\ca(\U(1))} & \stackrel{ \id\otimes\pi_{\cq_{\ca(\U(1))}} }{\longrightarrow} 
&\ASU\otimes (\ker\varepsilon_{\ca(\U(1))}/\mathcal{Q}_{\ca(\U(1))}) \,.
\end{array}
$$ 
Having a quantum homogeneous bundle, that is a quantum bundle whose total space is a Hopf algebra and whose fiber is a Hopf subalgebra of it, with the differential calculus on the fiber obtained from the corresponding projection, for the above sequence to be exact it is enough \cite{BM97} to check two conditions. The first one, 
$$
(\id\otimes\pi)\circ\Ad_{R}(\mathcal{Q}_{\SU}) \,\subset\, \mathcal{Q}_{\SU}\otimes \ca(\U(1)),
$$
with $\pi: \ASU \to \ca(\U(1))$ the projection in \eqref{qprp},
is easily established by a direct calculation and using the explicit form of the elements in $\mathcal{Q}_{\SU}$. The second condition amounts to the statement that the kernel of the projection $\pi$ can be written as a right $\ASU$-module of the kernel of $\pi$ itself restricted to the base algebra $\Asq$. Then, one needs to show that 
$
\ker\pi\,\subset\, (\ker\pi|_{\sq}) \ASU 
$,
the opposite implication being obvious. With $\pi$ defined in \eqref{qprp}, one has that 
$$
\ker\pi=\{ cf, \,c^{*}g, \quad \mathrm{with} \quad f,g \in\ASU\}. 
$$
Then $cf=c (a^{*}a+c^{*}c) f = ca^{*} (af) + c^{*}c (cf)$, with both $ca^{*}$ and $c^{*}c$ in $\ker\pi|_{\sq}$. The same holds for elements of the form $c^{*}g$, and the inclusion follows. 

We finish by showing that two of the generators in \eqref{q3dom} of the 3D calculus on $\ASU$, that is $\omega_\pm$, are indeed the generators of the horizontal forms, $\ker \sim_{\cn_{\SU}}$, on the principal bundle  
 as in \eqref{codd}. If $f\in\mathcal{A}$ with universal derivative $\delta f=1\otimes f-f\otimes 1$, from the definition \eqref{chimap} we get $ \chi(\delta f) = f_{(1)}\otimes\pi(f_{(2)})-f\otimes 1$, with the usual Sweedler notation:  $\Delta f=f_{(1)}\otimes f_{(2)}$. In particular, for the generators of $\ASU$ we get: 
\begin{align*}
\chi(\delta a)&=a\otimes(z-1), \quad 
\chi(\delta a^{*}) =a^{*}\otimes(z^{*}-1)
\nn \\
\chi(\delta c)&=c\otimes(z-1), \quad 
\chi(\delta c^{*})=c^{*}\otimes(z^{*}-1) .
\end{align*} 
Given the two generators $\omega_{\pm}$ and the specific $\cq_{\SU}$ which determines the 3D calculus, corresponding universal 1-forms can be taken to be:
$$
a\delta c-qc\delta a \in[\pi_{\cq_{\SU}}]^{-1}(\omega_{+}), \qquad 
c^* \delta a^*- qa^*\delta c^* \in[\pi_{\cq_{\SU}}]^{-1}(\omega_{-}).
$$
The action of the canonical map then gives:
\begin{align*}
&\chi(a\delta c-qc\delta a)=(ac-qca)\otimes(z-1)=0, \\ 
&\chi(c^* \delta a^*- qa^*\delta c^*) = (c^*a^* - qa^*c^*)\otimes(z^{*}-1)=0,
\end{align*}
which means that $\sim_{\cn_{\SU}}(\omega_{+})=0=\sim_{\cn_{\SU}}(\omega_{-})$. 
For the third generator $\omega_{z}$, one shows in a similar fashion that
$\sim_{\cn_{\SU}}(\omega_{z})=1\otimes(z-1)$. {}From these we may
conclude that the elements $\omega_\pm$ generate the module of horizontal forms. 

Finally, we know from \eqref{vvf} that the vector 
$X=X_z=(1-q^{-2})^{-1}(1-K^{4})$ is the dual generator to the calculus on the structure Hopf algebra $\ca(\U(1))$. For the corresponding `vector field' $\tilde{X}$ on $\ASU$ as in \eqref{hf}, one has that 
$\tilde{X}(\omega_\pm) = \langle{X},{\sim_{\cn_{\SU}}(\omega_{\pm})}\rangle=0$, while  
$\tilde{X}(\omega_z) = \langle{X},{\sim_{\cn_{\SU}}(\omega_{z})}\rangle=1$. 
These results identify $\tilde{X}$ as a vertical vector field.

%\newpage

\end{document}